\newif\ifpreprint
\preprinttrue 

\ifpreprint
	\documentclass[a4paper,12pt]{article}
	\usepackage{fullpage}
\fi

\usepackage{amsmath}
\usepackage{amssymb,amsfonts}
\ifpreprint
	\usepackage{amsthm}
\fi
\usepackage{enumitem}
\usepackage{graphicx}
\usepackage{xcolor}
\definecolor{darkblue}{RGB}{0,102,153}
\definecolor{cb2blue}{RGB}{55,126,184}
\definecolor{cb2green}{RGB}{77,175,74}
\definecolor{cb2red}{RGB}{228,26,28}
\definecolor{forestgreen}{RGB}{0,110,51}

\usepackage[colorlinks,
	citecolor=darkblue,
	linkcolor=darkblue,
	urlcolor=darkblue]{hyperref}

\usepackage[ruled,vlined,linesnumbered,algosection,algo2e]{algorithm2e} 

\usepackage[capitalize,noabbrev,nameinlink]{cleveref} 

\ifpreprint
	\theoremstyle{plain}
	\newtheorem{theorem}{Theorem}[section]
	\newtheorem{corollary}[theorem]{Corollary}
	\newtheorem{lemma}[theorem]{Lemma}
	\newtheorem{proposition}[theorem]{Proposition}
	\newtheorem{assum}{Assumption}
	\theoremstyle{definition}
	\newtheorem{definition}[theorem]{Definition}
	
	\newtheorem{remark}{Remark}[section]
\fi

\crefname{assum}{Assumption}{Assumptions}

\newcommand{\KeywordsAnd}{\and}
\ifpreprint
	\newcommand{\keywords}[1]{\par\noindent{\def\and{\unskip,\ }{\bf Keywords. }#1}\par}
	\newcommand{\subclass}[1]{\par\noindent{\def\and{\unskip,\ }{\bf AMS subject classifications. }#1}\par}
\fi
\ifpreprint
	\newcommand{\NL}[1][]{\ifstrempty{#1}{,}{#1}\ }
\fi

\crefname{line}{Step}{Steps} 
\SetKw{KwOr}{or}
\SetKw{KwAnd}{and}

\newcommand{\emailLink}[1]{\textsc{email} \href{mailto:#1}{#1}}
\newcommand{\orcidLink}[1]{\textsc{orcid} \href{https://orcid.org/#1}{#1}}
\newcommand{\amsmscLink}[1]{\href{http://www.ams.org/mathscinet/msc/msc2020.html?t=#1}{#1}}

\newcommand{\N}{\mathbb{N}}
\newcommand{\R}{\mathbb{R}}
\newcommand{\Rinf}{\overline{\R}}
\DeclareMathOperator*{\argmin}{\arg\min}
\DeclareMathOperator{\dom}{dom}
\DeclareMathOperator{\proj}{\Pi}
\DeclareMathOperator{\prox}{prox}
\DeclareMathOperator{\indicator}{\iota}
\DeclareMathOperator{\identity}{Id}
\newcommand{\coloneqq}{:=}
\DeclareMathOperator*{\minimize}{minimize}
\DeclareMathOperator{\stt}{subject~to}

\newcommand{\func}[3]{#1 \colon #2 \to #3}
\newcommand{\ffunc}[3]{#1 \colon #2 \rightrightarrows #3}
\newcommand{\innprod}[2]{\langle #1, #2 \rangle}
\newcommand{\normalcone}{\mathcal{N}}
\newcommand{\limnormalcone}{\normalcone^{\text{lim}}}
\DeclareMathOperator\lev{lev}
\DeclareMathOperator\dist{dist}
\DeclareMathOperator\cone{cone}
\DeclareMathOperator\epi{epi}
\DeclareMathOperator\interior{int}
\DeclareMathOperator\closedball{\mathbb{B}}
\newcommand{\jac}[1]{#1^\prime}

\ifpreprint
	\usepackage{tcolorbox}
	\newtcolorbox{mybox}[1][]{%
		left=0pt,
		right=0pt,
		top=0pt,
		bottom=0pt,
		colback=gray!12,
		colframe=gray!12,
		width=\dimexpr\textwidth\relax,
		enlarge left by=0mm,
		boxsep=4pt,
		arc=4pt,outer arc=4pt,
		#1
	}
\fi

\newcommand{\XX}{\R^n}
\newcommand{\YY}{\R^m}
\newcommand{\toattentive}[1]{\overset{#1}{\to}}
\newcommand{\LL}{\mathcal{L}}
\newcommand{\LLslack}{\LL^{\text{S}}}
\DeclareMathOperator{\Sproj}{\FBO}
\newcommand{\Ybounded}{Y}
\newcommand{\FBO}{\mathbf{O}}
\newcommand{\subLL}{p}
\newcommand{\subLLslack}{P}
\newcommand{\subLLsmooth}{F}
\newcommand{\subLLnm}{\Phi}
\newcommand{\Dvc}{D_{\text{VC}}}
\newcommand{\alps}{{ALPS}}
\newcommand{\als}{{ALS}}

\usepackage{pgfplots}
\pgfplotsset{compat=1.13}
\usepackage{tikz}
\usetikzlibrary{shapes,external}
\makeatletter
\providecommand{\leftsquigarrow}{%
	\mathrel{\mathpalette\reflect@squig\relax}%
}
\newcommand{\reflect@squig}[2]{%
	\reflectbox{$\m@th#1\rightsquigarrow$}%
}
\makeatother

\newcommand{\TheAuthorADM}{Alberto De~Marchi}
\newcommand{\TheEmailADM}{alberto.demarchi@unibw.de}
\newcommand{\TheOrcidADM}{0000-0002-3545-6898}
\newcommand{\TheAffiliationADM}{%
	University of the Bundeswehr Munich\NL
	Department of Aerospace Engineering\NL
	Institute of Applied Mathematics and Scientific Computing\NL
	85577 Neubiberg, Germany%
}
\newcommand{\TheAcknowledgementsADM}{%
	I wish to thank
	Andreas Themelis (Kyushu University) and
	Patrick Mehlitz (BTU Cottbus-Senftenberg),
	for critical and constructive comments on an early draft of this manuscript.
}
\newcommand{\TheTitle}{Implicit augmented Lagrangian and generalized optimization}

\newcommand{\TheKeywords}{%
	Nonsmooth nonconvex optimization\KeywordsAnd
	generalized nonlinear programming\KeywordsAnd
	multiplier methods\KeywordsAnd
	generalized augmented Lagrangians\KeywordsAnd
	proximal algorithms%
}
\newcommand{\TheAMSsubj}{%
	\amsmscLink{65K10}\and
	\amsmscLink{90C26}\and
	\amsmscLink{65K05}
}

\ifpreprint
	
\begin{document}
	
	\title{\bfseries \TheTitle}
	\author{\TheAuthorADM\thanks{%
			\TheAffiliationADM\NL[.]%
			\emailLink{\TheEmailADM}\NL%
			\orcidLink{\TheOrcidADM}\NL[.]%
		}%
	}%
	\date{\today}
\fi
	
\maketitle
\begin{abstract}
	Generalized nonlinear programming is considered without any convexity assumption, capturing a variety of problems that include nonsmooth objectives, combinatorial structures, and set-membership nonlinear constraints.
We extend the augmented Lagrangian framework to this broad problem class, preserving an implicit formulation and introducing slack variables merely as a formal device.
This, however, gives rise to a generalized augmented Lagrangian function that lacks regularity, due to the marginalization with respect to slack variables.
Based on parametric optimization, we develop a tailored stationarity concept to better qualify the iterates, generated as approximate solutions to a sequence of subproblems.
Using this variational characterization and the lifted representation, a suitable multiplier update rule is derived, and then asymptotic properties and convergence guarantees are established for a safeguarded augmented Lagrangian scheme.
An illustrative numerical example showcases the modelling versatility gained by dropping convexity assumptions and indicates the practical benefits of the advocated implicit approach.
\end{abstract}
\subclass{\TheAMSsubj}

\ifpreprint
	\clearpage
\fi

\section{Introduction}\label{sec:introduction}
Optimization problems often bring plenty of mathematical structure to work with.
In this paper we investigate the setting of generalized programming without any convexity, namely finite-dimensional optimization problems of the form
\begin{equation}
	\minimize_{x\in\XX}
	{}\quad{}
	\varphi(x) \coloneqq f(x) + g( c(x) ) ,
	\tag{P}\label{eq:P}
\end{equation}
where $x$ is the decision variable, $f$ and $c$ are smooth mappings and $g$ is merely lower semicontinuous.
Problem data $f$ and $g$ are allowed to be nonconvex mappings, the nonsmooth cost $g$ is not necessarily continuous nor real-valued, and the mapping $c$ is potentially nonlinear.
The composition of a nonsmooth term $g$ with a smooth mapping $c$ results in great modeling flexibility and encompasses a broad spectrum of problems \cite{rockafellar2000extended}; both regularizers and constraints can be encoded, as well as combinations thereof \cite{combettes2011proximal,boyd2011distributed,demarchi2021dissertation}.
Prominent special cases of \eqref{eq:P} are nonlinear \cite{bertsekas1996constrained} and disjunctive \cite{balas2018disjunctive} programming, structured optimization \cite{parikh2014proximal,themelis2018forward}, and constrained structured optimization \cite{demarchi2023constrained}.
The abstract problem \eqref{eq:P} is also referred to as ``extended nonlinear programming'' in \cite{rockafellar2000extended}, promoting better representation of structures found in applications by featuring a composite term, beyond the conventional format of nonlinear programming.

A fundamental technique for solving \emph{constrained} optimization problems is the augmented Lagrangian (AL) framework \cite{rockafellar1976augmented,bertsekas1996constrained,birgin2014practical}.
Introducing an auxiliary variable $z \in \YY$, \eqref{eq:P} can be equivalently rewritten in the form
\begin{equation}
	\minimize_{x\in\XX ,\ z\in\YY}
	{}\quad{}
	f(x) + g(z)
	{}\qquad{}
	\stt
	{}\quad{}
	c(x) - z = 0 ,
	\tag{P$_{\text{S}}$}\label{eq:ALMz:P}
\end{equation}
which has a simpler, separable objective function but is subject to some (explicit) constraints.
It has been observed in \cite[Lemma 3.1]{demarchi2023constrained} that incorporating slack variables in this manner does not affect minimizers and stationary points when compared to the original \eqref{eq:P} --- somewhat remarkably in light of \cite{benko2021implicit}.
Constrained minimization problems such as \eqref{eq:ALMz:P} are amenable to be addressed by means of AL methods, which we embrace as they can effortlessly handle nonsmoothness \cite{dhingra2019proximal,evens2021neural}.
As the AL approach constitutes a framework, rather than a single algorithm, several methods have been presented in the past decades, expressing the foundational ideas in different flavors.

Since their inception \cite{hestenes1969multiplier,powell1969method,rockafellar1973dual}, AL schemes have been deeply and extensively investigated \cite{conn1991globally,bertsekas1996constrained,gill2012primal,birgin2014practical,grapiglia2020complexity}, also in infinite dimensions \cite{kanzow2018augmented}.
In the convex setting, it was recognized that the ``method of multipliers'' \cite{hestenes1969multiplier,powell1969method} corresponds to solving an associated dual problem via the proximal point algorithm (PPA) \cite{rockafellar1973dual,rockafellar1976augmented}, with duality theory playing a key role in this parallel.
Then, sufficient optimality conditions and local duality led to analogous results in nonlinear programming by Bertsekas \cite{bertsekas1996constrained}.
Pivoting on some kind of local duality, in turn enabled by some local convexity, the AL scheme can be traced back to an application of the PPA, thereby establishing convergence.

This pattern was recently followed by Rockafellar \cite{rockafellar2022augmented,rockafellar2022convergence} to extend these ideas to the much broader setting of \eqref{eq:P} with $g$ convex.
In contrast with this line of work, which goes back at least to
\cite[Section 2]{rockafellar2000extended}, we are interested here in investigating AL methods for generalized nonlinear programming without any convexity assumption.
With $g$ potentially nonconvex, the abstract setting \eqref{eq:P} not only covers a variety of problems, but it also provides a more versatile and expressive modeling paradigm, with no need for ``suitable reformulations''.
Dropping the convexity of $g$ comes at a price, though, as it appears more difficult to leverage the perspective of PPA.
Nonetheless, the shifted-penalty approach underpinning the seminal ``method of multipliers'' still applies, and can be exploited as such.
Yet difficulties arise as the nonconvexity of $g$ translates into potential set-valuedness of its proximal mapping, which in turn leads to a lack of regularity.
Here, we seek a better understanding of the consequences on the analysis and design of numerical methods for \eqref{eq:P}, concerning theoretical properties, convergence guarantees, and practical aspects.

The following blanket assumptions are considered throughout, without further mention.
Technical definitions are given in \cref{sec:preliminaries}.
\begin{mybox}
	\begin{assum}\label{ass:P}
		The following hold in \eqref{eq:P}:
		\begin{enumerate}[label=(\roman*)]
			\item\label{ass:f}%
			$\func{f}{\XX}{\R}$ and $\func{c}{\XX}{\YY}$ are continuously differentiable;
			\item\label{ass:g}%
			\(\func{g}{\YY}{\Rinf}\) is proper, lower semicontinuous, and prox-bounded;
			\item\label{ass:phi}%
			\(\inf\varphi\in\R\).
		\end{enumerate}
	\end{assum}
\end{mybox}
\cref{ass:P}\ref{ass:phi} is relevant to the well-posedness of the problem, as the cost function $\varphi \coloneqq f + g \circ c$ has nonempty domain and is bounded from below.
Prox-boundedness in \cref{ass:P}\ref{ass:g} guarantees that the proximal mapping of $g$ is well-defined, for some suitable stepsizes.
Moreover, we work under the practical assumption that (only) the following computational oracles are available or simple to evaluate:
\begin{itemize}
	\item given $x \in \XX$, cost function value $f(x) \in \R$ and gradient $\nabla f(x) \in \XX$;
	\item given $v \in \YY$ and $\gamma \in (0,\gamma_g)$, with $\gamma_g$ being the prox-boundedness threshold of $g$, proximal point $z \in \prox_{\gamma g}(v)$ and cost function value $g(z) \in \R$ therein;
	\item given $x \in \XX$ and $v \in \YY$, constraint function value $c(x) \in \YY$ and transposed Jacobian-vector product $\jac c(x)^\top v \in \XX$,
	where $\func{\jac c}{\XX}{\R^{m \times n}}$ denotes the Jacobian of $c$.
\end{itemize}
Relying only on these primitives, the numerical method presented in the following is first-order and matrix-free by construction; as such, it involves only simple operations and has low memory footprint.

\paragraph{Outline}
Related works and topics are mentioned in \cref{sec:relatedwork}, followed by our contributions in \cref{sec:contributions}.
Notation, preliminary material, and useful tools are given in \cref{sec:preliminaries}, before introducing relevant optimality and stationarity notions in \cref{sec:concepts}.
The proposed augmented Lagrangian scheme is detailed and characterized in \cref{sec:ALM}.
A solver for the inner AL problems is developed and analysed in \cref{sec:inner}.
Finally, an illustrative numerical example is reported in \cref{sec:numresults}, before some concluding remarks in \cref{sec:conclusions}.

	\subsection{Related Work}\label{sec:relatedwork}
	The template \eqref{eq:P} covers many classical problem classes.
These include structured two-term problems $f + g$ \cite{parikh2014proximal} and
three-term problems $f + g + h$, as well as settings with the nonsmooth term composed with a linear operator, with two-term $g + h \circ L$, three-term problems $f + g + h \circ L$ \cite{latafat2017asymmetric}, or with linear constraints as for the alternating direction method of multipliers \cite{boyd2011distributed}.
Most of these problem classes have been extensively studied in the convex setting, while nonconvexity has been mostly considered for structured two-term problems $f + g$ only \cite{themelis2018forward}.

In conjunction with AL schemes, non-Lipschitzian or merely lower semicontinuous cost functions have been recently considered in \cite{chen2017augmented,evens2021neural,demarchi2023constrained}.
The issue of attentive convergence, defined in \cref{sec:preliminaries} and seemingly necessary to perform limiting analysis, was overlooked in \cite{chen2017augmented} and brought to attention in \cite{demarchi2023constrained}.
Since it is difficult to attest a priori, as a sufficient to guarantee attentive convergence one may additionally assume $g$ continuous relative to its domain.
While much weaker than convexity, this assumption is satisfied by many practical problems, as pointed out in \cite{demarchi2022interior}.

Another issue is that of well-definedness of the AL algorithm.
In fact, the AL subroblems could be ill-posed, as there is no a priori guarantee that the AL function is bounded from below, due to the constraints relaxation \cite{birgin2014practical}.
As we are concerned with the role of $g$, in this work, we neglect this possibility and assume that all subproblems are well-posed.
Several approaches for moving beyond this assumption may be applicable, although possibly imposing other requirements.
Among these are localized minimizations by exploiting non-relaxable constraints \cite{andreani2008augmented} or a trust-region included in the subproblems, as well as different AL schemes \cite{grapiglia2020complexity,evens2021neural}.

Constrained structured optimization \cite{demarchi2023constrained} entails the minimization of the sum $f(x) + g(x)$ over $x \in \XX$ subject to set-membership constraints $c(x) \in D$, where $c$ is a smooth mapping and $D$ a nonempty closed set.
As such, it is more general than nonlinear and disjunctive programming \cite{bertsekas1996constrained,balas2018disjunctive}, but covered by the template \eqref{eq:P}.
Conversely, possibly introducing some auxiliary variables, \eqref{eq:P} can always be reformulated as a constrained structured program.
Thus, the problem class investigated in \cite{demarchi2023constrained} is as general as \eqref{eq:P} under \cref{ass:P}, but our solution strategy here is essentially different.
	
	\subsection{Approach and Contributions}\label{sec:contributions}
	The core idea of this work traces back to the author's doctoral dissertation \cite[Chapter 1]{demarchi2021dissertation}.
With the aim of improving the analysis and patching weak reasoning, especially resolving the issue of a suitable termination criterion for the subsolver, further developments took a different spirit and eventually became \cite{demarchi2023constrained}.
Therein slack variables are introduced as in \eqref{eq:ALMz:P}, making the constraints \emph{explicit} and then applying a traditional AL scheme.
To the contrary, here we maintain the \emph{implicit} approach advocated for in \cite{demarchi2021dissertation}, offering a rigorous convergence analysis.
In particular, we use the lifted reformulation \eqref{eq:ALMz:P} only as a theoretical device to develop our strategy for solving \eqref{eq:P} and investigate its properties via sharpened variational tools.

Marginalizing the AL function of \eqref{eq:ALMz:P} with respect to the slack variable $z$ yields the generalized (or, implicit) AL function for \eqref{eq:P} \cite{dhingra2019proximal,demarchi2021dissertation,rockafellar2022convergence}.
However, by the nonconvexity of $g$,
the generalized AL function fails to be continuously differentiable (or prox-friendly) in general.
For overcoming this lack of regularity,
\begin{itemize}
	\item
	we devise the concept of $\Upsilon$-stationarity, which allows us to characterize precisely the sought (approximate) solutions to the AL subproblems,
\end{itemize}
thereby better capturing the intended philosophy.
Rooted in parametric optimization, $\Upsilon$-stationarity is strong enough for the outer AL scheme
but still computable and easy to handle numerically.
Indeed, we show this concept provides an intermediate qualification between M-stationarity for the implicit and explicit formulations.

With the concept of $\Upsilon$-stationarity at hand, we focus on
solving the AL subproblems, namely minimizing the AL function, which constitutes the major computational task and a key step in AL algorithms.
We explore the possibility to directly solve the AL subproblems, without slack variables.
In particular,
\begin{itemize}
	\item
	we design and analyze a subsolver that exploits the proximal oracle of $g$ for minimizing the generalized AL function and, by generating approximate $\Upsilon$-stationary points, complies with the requirements of the AL scheme.
\end{itemize}
Reminiscent of \emph{magical steps} \cite[Section 10.4.1]{conn2000trust}, the rationale behind the subsolver is to include slack variables only as algorithmic by-products, as opposed to decision variables.
These are then adopted within an iterative scheme to generate updates that yield sufficient improvement and to terminate with a suitable certificate of (approximate) stationarity.

Overall, the advantages of weaker assumptions and stronger stationarity concepts would not be attractive without convergence results.
Therefore,
\begin{itemize}
	\item
	we investigate the asymptotic properties of a safeguarded implicit AL scheme for \eqref{eq:P} under \cref{ass:P} and establish (global) convergence guarantees on par with those in classical nonlinear programming.
\end{itemize}
Moreover, we illustrate how following an implicit approach without introducing slack variables can result in improved modeling flexibility and practical performance.
Not only does sticking to a compact problem formulation avoid increasing the problem size, but it also leads to better use of the available computational oracles and tends to avoid spurious (local) minima and stationary points.

\section{Notation and Fundamentals}\label{sec:settingfundamentals}
In this section, we comment on notation, preliminary definitions and useful results.

	\subsection{Preliminaries}\label{sec:preliminaries}
	With $\R$ and $\Rinf \coloneqq \R \cup \{\infty\}$ we denote the real and extended-real line, respectively.
The vectors in $\R^n$ with all elements equal to $0$ or $1$ are denoted as $0_n$ and $1_n$; whenever $n$ is clear from context we simply write $0$ and $1$, respectively.
With $\closedball_r(x)$ we indicate the closed ball centered at $x$ with radius $r$.
The \emph{effective domain} of an extended-real-valued function $\func{h}{\R^n}{\Rinf}$ is denoted by $\dom h \coloneqq \{ x\in\R^n \,\vert\, h(x) < \infty \}$.
We say that $h$ is \emph{proper} if $\dom h \neq \emptyset$
and \emph{lower semicontinuous} (lsc) if $h(\bar{x}) \leq \liminf_{x\to\bar{x}} h(x)$ for all $\bar{x} \in \R^n$.
For some constant $\tau\in\R$, $\lev_{\leq \tau} h \coloneqq \{ x\in\XX \,\vert\, h(x)\leq \tau \}$
denotes the $\tau$-\emph{sublevel set} associated with $h$.
Function $h$ is (lower) \emph{level-bounded} if for every $\tau \in \R$ its $\tau$-sublevel set is bounded.
Then, we introduce the following notion of \emph{uniform level boundedness}, which plays an important role in parametric minimization \cite[Definition 1.16]{rockafellar1998variational}.
A function $\func{f}{\XX \times \YY}{\Rinf}$ with values $f(x,u)$ is \emph{level-bounded in
$x$ locally uniformly in $u$} if for each $\alpha \in \R$ and $\bar{u} \in \YY$ there exists $\varepsilon > 0$ such that the set $\{ (x,u) \,\vert\, f(x,u) \leq \alpha, \|u - \bar{u}\| \leq \varepsilon \}$ is bounded in $\XX \times \YY$.
We use the notation $\ffunc{F}{\XX}{\YY}$ to indicate a point-to-set function $\func{F}{\XX}{\mathcal{P}(\YY)}$, where $\mathcal{P}(\YY)$ is the power set of $\YY$ (the set of all subsets of $\YY$).

Given a proper and lsc function $\func{h}{\XX}{\Rinf}$ and a point $\bar{x}$ with $h(\bar{x})$ finite, we may avoid to assume $h$ continuous and instead appeal to $h$-\emph{attentive} convergence of a sequence $\{ x^k \}$, denote as $x^k \toattentive{h} \bar{x}$ and given by $x^k \to \bar{x}$ with $h(x^k) \to h(\bar{x})$.
Following \cite[Definition 8.3]{rockafellar1998variational}, we denote by $\ffunc{\hat{\partial} h}{\XX}{\XX}$ the \emph{regular subdifferential} of $h$, where
\begin{equation*}
	v \in \hat{\partial} h(\bar{x})
	{}\quad:\Leftrightarrow\quad{}
	\liminf_{\substack{x\to\bar{x}\\x\neq\bar{x}}} \frac{h(x) - h(\bar{x}) - \langle v, x-\bar{x}\rangle}{\|x-\bar{x}\|} \geq 0 .
\end{equation*}
The (limiting) \emph{subdifferential} of $h$ is $\ffunc{\partial h}{\XX}{\XX}$, where $v \in \partial h(\bar{x})$ if and only if $\bar{x} \in \dom h$ and there exist sequences $\{x^k\}$ and $\{v^k\}$ such that $x^k \toattentive{h} \bar{x}$ and $v^k \in \hat{\partial} h(x^k)$ with $v^k \to v$.
The subdifferential of $h$ at $\bar{x}$ satisfies $\partial(h+h_0)(\bar{x}) = \partial h(\bar{x}) + \nabla h_0(\bar{x})$ for any $\func{h_0}{\XX}{\Rinf}$ continuously differentiable around $\bar{x}$ \cite[Exercise 8.8]{rockafellar1998variational}.

Given a parameter value $\gamma>0$, the \emph{proximal} mapping $\prox_{\gamma h}$ is defined by
\begin{equation}
	\prox_{\gamma h}(x)
	{}\coloneqq{} 
	\argmin_{z} \left\{ h(z) + \frac{1}{2\gamma}\|z-x\|^2 \right\} .
	\label{eq:proxmapping}
\end{equation}
We say that $h$ is \emph{prox-bounded} if it is proper and $h + \|\cdot\|^2 / (2\gamma)$ is bounded below on $\XX$ for some $\gamma > 0$ \cite[Definition 1.23]{rockafellar1998variational}.
The supremum of all such $\gamma$ is the threshold $\gamma_h$ of prox-boundedness for $h$.
In particular, if $h$ is bounded below by an affine function, then $\gamma_h = \infty$.
When $h$ is lsc, for any $\gamma \in (0,\gamma_h)$ the proximal mapping $\prox_{\gamma h}$ is locally bounded, nonempty- and compact-valued \cite[Theorem 1.25]{rockafellar1998variational}.
The value function of the minimization problem defining the proximal mapping is the \emph{Moreau envelope} with parameter $\gamma > 0$, denoted $\func{h^\gamma}{\XX}{\R}$, namely
\begin{equation*}
	h^\gamma(x) \coloneqq \inf_{z} \left\{ h(z) + \frac{1}{2\gamma}\|z-x\|^2 \right\} .
\end{equation*}

Some tools of variational analysis will be exploited in order to describe the geometry of some nonempty, closed, but not necessarily convex, set $D \subset \YY$.
The \emph{projection} mapping $\ffunc{\proj_D}{\YY}{\YY}$ and the \emph{distance} function $\func{\dist_D}{\YY}{\R}$ are defined by
\begin{equation*}
	\proj_D(v)
	{}\coloneqq{}
	\argmin_{z \in D} \|z - v\|
	{}\qquad\text{and}\qquad{}
	\dist_D(v)
	{}\equiv{}
	\dist(v, D)
	{}\coloneqq{}
	\inf_{z \in D} \|z - v\| .
\end{equation*}
The former is a set-valued mapping whenever $D$ is nonconvex, whereas the latter is always single-valued.
The \emph{indicator} of $D$ is denoted $\func{\indicator_D}{\YY}{\Rinf}$ and defined by $\indicator_D(v) = 0$ if $v \in D$, and $\indicator_D(v) = \infty$ otherwise.
If $D$ is nonempty and closed, then $\indicator_D$ is proper and lsc.
The proximal mapping of $\indicator_D$ is the projection $\proj_D$; thus, $\proj_D$ is locally bounded.
Given $z \in D$, the \emph{limiting normal cone} to $D$ at $z$ is the closed cone
\begin{equation*}
	\limnormalcone_D(z) \coloneqq \limsup_{v \to z} \; \cone \left( v - \proj_D(v) \right) .
\end{equation*}
The limiting normal cone is robust in the sense that for all $z \in D$ it satisfies $\limnormalcone_D(z) = \limsup_{v \to z} \limnormalcone_D(v)$.
For any proper and lsc function $\func{h}{\XX}{\Rinf}$ and point $\bar{x} \in \dom h$, we have
\begin{equation}
	\label{eq:subdiffLimNormalConeEpigraph}
	\partial h(\bar{x})
	=
	\left\{v\in\XX \,\vert\, (v,-1)\in\limnormalcone_{\epi h}(\bar{x},h(\bar{x}))\right\} ,
\end{equation}
where $\epi h \coloneqq \{(x,\alpha) \in \XX \times \R \,\vert\, h(x) \leq \alpha \}$ denotes the \emph{epigraph} of $h$ \cite[Theorem 8.9]{rockafellar1998variational}.
	
	\subsection{Stationarity Concepts}\label{sec:concepts}
	We now define some basic concepts and discuss optimality and stationarity conditions for \eqref{eq:P}.
As the cost function $\varphi$ is possibly extended-real-valued, the concept of feasibility for a point refers to the domain of this cost function.
A point $x^\ast \in \XX$ is called \emph{feasible} for \eqref{eq:P} if $x^\ast \in \dom \varphi$, that is, if $c(x^\ast) \in \dom g$.
With \eqref{eq:P}, we seek points where the minimum of $\varphi$ is attained, which are necessarily feasible.
Notions of optimality can have global or local scope.
A feasible point $x^\ast \in \XX$ is said \emph{locally optimal}, or called a \emph{local minimizer}, for \eqref{eq:P} if there exists $\delta > 0$ such that $\varphi(x^\ast) \leq \varphi(x)$ holds for all feasible $x \in \closedball_\delta(x^\ast)$.
Additionally, if the inequality holds for all feasible $x \in \XX$, then $x^\ast$ is said \emph{(globally) optimal} for \eqref{eq:P}.

Focusing on the development of an ``affordable'' algorithm, our method is designed to seek \emph{candidate} minimizers, guided by necessary optimality conditions.
Working without any convexity assumptions, a plausible stationarity concept is that based on limiting subdifferentials; cf. \cite[Section 3]{mehlitz2020asymptotic} and \cite[Theorem 5.48]{mordukhovich2006variational}.
With respect to the minimization of a proper lsc function $\func{h}{\XX}{\Rinf}$, we say that $x^\ast \in \dom h$ is \emph{M-stationary} if $0 \in \partial h(x^\ast)$, which constitutes a necessary condition for the optimality of $x^\ast$ \cite[Theorem 10.1]{rockafellar1998variational}.
When addressing \eqref{eq:P}, the composite structure gives rise to first-order optimality conditions expressed in Lagrangian terms \cite{rockafellar2022augmented,rockafellar2022convergence}, as it intuitively follows from \eqref{eq:ALMz:P}. 
Although closely related to M-stationarity,
we prefer the nomenclature \emph{KKT}-stationarity
to emphasize the Lagrangian structure.
\begin{mybox}
	\begin{definition}[KKT-stationarity]\label{def:Mstationary}
		Relative to \eqref{eq:P}, a point $x^\ast \in \XX$ is called \emph{KKT-stationary} if there exists a multiplier $y^\ast \in \YY$ such that
		\begin{subequations}
			\begin{align}
				0
				{}={}&
				\nabla f(x^\ast) + \jac c(x^\ast)^\top y^\ast ,
				\label{eq:Mstationary:x} \\
				y^\ast
				{}\in{}&
				\partial g( c(x^\ast) ) .
				\label{eq:Mstationary:y}
			\end{align}
		\end{subequations}
	\end{definition}
\end{mybox}
Notice that \eqref{eq:Mstationary:y} implicitly requires the feasibility of $x^\ast$, namely $c(x^\ast) \in \dom g$, for otherwise the subdifferential $\partial g(c(x^\ast))$ is empty.
Moreover, note that \cref{def:Mstationary} coincides with the celebrated Karush-Kuhn-Tucker (KKT) conditions in nonlinear programming if $g$ in \eqref{eq:P} is the indicator of a nonempty closed convex set \cite{bertsekas1996constrained,rockafellar2022convergence,demarchi2023constrained}.

Although potentially ambiguous in the nonconvex setting, one could refer to the Lagrange multiplier $y$ as to the \emph{dual} variable and to $(x,y)$ as to a \emph{primal-dual} pair.
In fact, a sort of local duality emerged while investigating convergence of AL methods in nonconvex settings, as revealed by Bertsekas \cite{bertsekas1996constrained} and recently by Rockafellar in a broader context \cite{rockafellar2022convergence,rockafellar2022augmented}---but this seems premature for \eqref{eq:P} under \cref{ass:P}.

Given some tolerance $\varepsilon \geq 0$, an approximate M-stationarity concept for the minimization of $h$ refers to $\dist( 0, \partial h(x^\ast) ) \leq \varepsilon$.
We may also consider relaxed conditions for an approximate counterpart of \cref{def:Mstationary}.
\begin{mybox}
	\begin{definition}[Approximate KKT-stationarity]\label{def:approxMstationary}
		Let $\varepsilon \geq 0$ be given.
		Relative to \eqref{eq:P}, a point $x^\ast \in \XX$ is called \emph{$\varepsilon$-KKT-stationary} if there exist $z^\ast \in \dom g$ and a multiplier $y^\ast \in \YY$ such that
		\begin{subequations}\label{eq:approxMstationary}
			\begin{align}
				\| \nabla f(x^\ast) + \jac c(x^\ast)^\top y^\ast \|
				{}\leq{}&
				\varepsilon ,
				\label{eq:approxMstationary:x} \\
				\dist( y^\ast, \partial g( z^\ast ) )
				{}\leq{}&
				\varepsilon ,
				\label{eq:approxMstationary:y} \\
				\| c(x^\ast) - z^\ast \|
				{}\leq{}&
				\varepsilon .
				\label{eq:approxMstationary:z}
			\end{align}%
		\end{subequations}%
	\end{definition}
\end{mybox}
It should be noticed that $\varepsilon$-KKT-stationarity with $\varepsilon = 0$ recovers the notion of KKT-stationarity, analogously to M-stationarity.
This fact follows from $\partial g(z)$ being closed at any point $z \in \dom g$ \cite[Theorem 8.6]{rockafellar1998variational}, which implies that, for any $v \in \YY$, $\dist( v, \partial g(z) ) = 0$ is equivalent to the inclusion $v \in \partial g(z)$.

The approximation allowed by \cref{def:approxMstationary} entails not only stationarity \eqref{eq:approxMstationary:x} and complementarity \eqref{eq:approxMstationary:y}, but affects also feasibility: owing to \eqref{eq:approxMstationary:z} it may be that $c(x^\ast) \notin \dom g$.
Encompassing this inexactness is particularly important in practice as it might be difficult to deal with the composition $g \circ c$ as soon as $c$ is a nonlinear mapping.
On the contrary, it is dramatically easier to deal with \eqref{eq:P} when $c$ is the identity mapping.
In this particular case, we apply a narrower definition compared to \cref{def:approxMstationary}, as we deem feasibility of an approximate M-stationary point a significant property.
Therefore, relative to \eqref{eq:P} with $c \coloneqq \identity$, a point $x^\ast \in \XX$ is called \emph{$\varepsilon$-M-stationary} if $\dist( - \nabla f(x^\ast), \partial g( x^\ast ) ) \leq \varepsilon$.
Along with feasibility, namely $x^\ast \in \dom g$, this demands a tighter approximation than \eqref{eq:approxMstationary}.
As it suffices to invoke the proximal mapping of $g$ to reconcile the ``constraint'' $x\in\dom g$, many feasible methods exist for such two-term problems $f+g$; see e.g. \cite{combettes2011proximal,parikh2014proximal,themelis2018forward,demarchi2022proximal} and references therein.

We add to our theoretical toolbox also an asymptotic counterpart of KKT-stationarity.
The following closely patterns the \emph{approximate} KKT conditions of \cite[Definition 3.1]{birgin2014practical}, and so does the meaning of primal iterates $x^k$ and Lagrange multipliers $y^k$.
\begin{mybox}
	\begin{definition}[AKKT-stationarity]\label{def:AMstationary}
		Relative to \eqref{eq:P}, a feasible point $x^\ast \in \XX$ is called \emph{asymptotically KKT-stationary} if there exist sequences $\{x^k\}_{k\in\N} \subseteq \XX$ and $\{y^k\}_{k\in\N}, \{z^k\}_{k\in\N} \subseteq \YY$ such that $z^k \toattentive{g} c(x^\ast)$ and
		\begin{subequations}\label{eq:AMstationary}%
			\begin{align}%
				\nabla f(x^k) + \jac c(x^k)^\top y^k
				{}\to{}&
				0 ,
				\label{eq:AMstationary:x} \\
				y^k
				{}\in{}&
				\partial g( z^k ) ,
				\label{eq:AMstationary:y} \\
				c(x)^k - z^k
				{}\to{}&
				0 .
				\label{eq:AMstationary:z}
			\end{align}%
		\end{subequations}%
	\end{definition}
\end{mybox}
Notice that \cref{def:AMstationary} does not require the sequence $\{ y^k \}_{k\in\N}$ to converge.
A measure of (primal) infeasibility is suggested by \eqref{eq:AMstationary:z}, where $z^k$ is introduced to account in \eqref{eq:AMstationary:y} for the fact that the condition $c(x^k) \in \dom g$ can be violated along the iterates, though it should hold asymptotically.
Similarly, with \eqref{eq:AMstationary:x} one can monitor the inexactness in the stationarity condition, or \emph{dual} infeasibility.
Finally, the attentive convergence $z^k \toattentive{g} c(x^\ast)$ plays a role when taking the limit in \eqref{eq:AMstationary}, since we aim to recover the limiting subdifferential for KKT-stationarity; cf. \cref{def:Mstationary}.

A local minimizer for \eqref{eq:P} is KKT-stationary under validity of a suitable qualification condition, which, by non-Lipschitzness of $g$, will depend on function $g$ as well, see \cite{guo2018necessary}.
However, we can show that each local minimizer of \eqref{eq:P} is always AKKT-stationary.
Related results can be found in \cite[Theorem 6.2]{kruger2021optimality}, \cite[Section 5.1]{mehlitz2020asymptotic} and \cite[Proposition 2.5]{demarchi2023constrained}.
\begin{mybox}
	\begin{proposition}\label{prop:minimizers_AM_stationary}
		Let $x^\ast \in\XX$ be a local minimizer for \eqref{eq:P}.
		Then, $x^\ast$ is AKKT-stationary.
	\end{proposition}
\end{mybox}
\begin{proof}\label{prop:minimizers_AM_stationary:proof}
	By local optimality of $x^\ast$ for \eqref{eq:P} there exists $\delta > 0$ such that $\varphi(x^\ast) \leq \varphi(x)$ is valid for all feasible $x \in \closedball_\delta(x^\ast)$.
	Consequently, $x^\ast$ is the unique global minimizer of the localized problem
	\begin{equation}\label{eq:regularized_problem}
		\minimize_{x \in \closedball_\delta(x^\ast)}
		{}\quad{}
		\varphi(x) + \frac{1}{2}\|x-x^\ast\|^2 .
	\end{equation}
	Let us now consider the penalized surrogate problem
	\begin{equation}\label{eq:penalized_surrogate}\tag{P$(k)$}
		\minimize_{\substack{x \in \closedball_\delta(x^\ast)\\z \in \closedball_1(c(x^\ast))}}
		{}\quad{}
		f(x) + g(z) + \frac{\varrho_k}{2}\|c(x) - z\|^2 + \frac{1}{2}\|x-x^\ast\|^2
	\end{equation}
	where $k\in\N$ is arbitrary, $\varrho_k >0$, and the sequence $\{\varrho_k\}_{k\in\N}$ satisfies $\varrho_k\to\infty$ as $k\to\infty$.
	Noting that the objective function of this optimization problem is lsc while its feasible set is nonempty (by feasibility of $x^\ast$) and compact, it possesses a global minimizer $(x^k, z^k)$ for each $k\in\N$.
	Therefore, without loss of generality, up to extracting a subsequence, we assume $x^k\to\widehat{x}$ and $z^k\to\widehat{z}$ for some $\widehat{x}\in\closedball_\delta(x^\ast)$ and $\widehat{z}\in \dom g \cap \closedball_1(c(x^\ast))$.

	We now argue that $\widehat{x}=x^\ast$ and $\widehat{z}=c(\widehat{x})=c(x^\ast)$.
	To this end, we note that $(x^\ast,c(x^\ast))$ is feasible to \eqref{eq:penalized_surrogate},
	which yields the estimate
	\begin{equation}\label{eq:estimate_from_penalization}
		f(x^k) + g(z^k) + \frac{\varrho_k}{2}\|c(x^k) - z^k\|^2 + \frac12\|x^k-x^\ast\|^2
		\leq
		\varphi(x^\ast)
	\end{equation}
	for each $k\in\N$.
	Using lower semicontinuity of $g$, finiteness of $\varphi(x^\ast)$ as well as the convergences $c(x^k)\to c(\widehat{x})$ and $z^k\to\widehat{z}$, taking the limit for $k\to\infty$ in \eqref{eq:estimate_from_penalization} gives $c(\widehat{x})=\widehat{z}\in \dom g$.
	Therefore, $\widehat{x}$ is feasible for \eqref{eq:regularized_problem}
	and local optimality of $x^\ast$ for \eqref{eq:P} implies $\varphi(x^\ast) \leq \varphi(\widehat{x})$.
	Furthermore, exploiting \eqref{eq:estimate_from_penalization} and the optimality of each $(x^k,z^k) \in \XX \times \dom g$, we find
	that $\widehat{x}=x^\ast$, since
	\begin{multline*}
		\varphi(\widehat{x})+\frac12\|\widehat{x}-x^\ast\|^2
		\leq
		\liminf\limits_{k\to\infty}
		\left(
		f(x^k) + g(z^k) + \frac{\varrho_k}{2}\|c(x^k)-z^k\|^2+\frac12\|x^k-x^\ast\|^2
		\right) \\
		\leq
		\varphi(x^\ast)
		\leq
		\varphi(\widehat{x}).
	\end{multline*}
	Moreover, noting that \eqref{eq:estimate_from_penalization} also gives $f(x^k) + g(z^k) \leq \varphi(x^\ast)$ for each $k\in\N$, by
	\begin{align*}
		\varphi(x^\ast)
		\leq
		\liminf\limits_{k\to\infty} f(x^k) + g(z^k)
		\leq
		\limsup\limits_{k\to\infty} f(x^k) + g(z^k)
		\leq
		\varphi(x^\ast)
	\end{align*}
	we have the attentive convergence of $\{z^k\}_{k\in\N}$ to $\widehat{z} = c(x^\ast)$,
	that is, $z^k \toattentive{g} c(x^\ast)$.
	
	Now, due to $x^k\to x^\ast$ and $z^k\to c(x^\ast)$, we may assume without loss of generality that $\{x^k\}_{k\in\N}$ and $\{z^k\}_{k\in\N}$ are taken from the interior of $\closedball_\delta(x^\ast)$ and $\closedball_1(c(x^\ast))$, respectively, as this is eventually the case.
	Thus, for each $k\in\N$, $(x^k, z^k)$ is an unconstrained (that is, with the constraints being inactive) minimizer of
	\[
	(x, z)\mapsto f(x) + g(z) + \frac{\varrho_k}{2}\|c(x) - z\|^2 + \frac{1}{2}\|x - x^\ast\|^2 ,
	\]
	whose relevant necessary conditions for optimality read
	\begin{align*}
		0 {}={}& \nabla f(x^k) + \varrho_k \jac c(x^k)^\top [ c(x^k) - z^k ] + x^k - x^\ast , \\
		0 {}\in{}& \partial g(z^k) - \varrho_k [ c(x^k) - z^k ] .
	\end{align*}
	Thus, setting $y^k \coloneqq \varrho_k [ c(x^k) - z^k ]$ for each $k\in\N$, the conditions in \eqref{eq:AMstationary} are a consequence of $x^k \to x^\ast$ and $z^k \to c(x^\ast)$.
	Overall, this shows that $x^\ast$ is AKKT-stationary for \eqref{eq:P}.
	\qedhere
\end{proof}
In order to guarantee that local minimizers for \eqref{eq:P} are not only AKKT- but also KKT-stationary, the presence of a qualification condition is necessary.
Analogously to \cite[Definition 2.6]{demarchi2023constrained}, the following regularity notion generalizes the (comparatively weak) constraint qualification from \cite[Section 3.2]{mehlitz2020asymptotic} to the non-Lipschitzian setting of \eqref{eq:P} and is closely related to the \emph{uniform qualification condition} introduced in \cite[Definition 6.8]{kruger2021optimality}.
\begin{mybox}
	\begin{definition}[AKKT-regularity]\label{def:AMregularity}
		Let $\bar{x} \in \XX$ be a feasible point for \eqref{eq:P},
		namely $c(\bar{x}) \in \dom g$,
		and define the set-valued mapping $\ffunc{\mathcal{M}}{\XX \times \YY}{\XX}$ by $\mathcal{M}(x,z) \coloneqq \jac c(x)^\top \partial g(z)$.
		Then, $\bar{x}$ is called \emph{asymptotically KKT-regular} for \eqref{eq:P} if
		\begin{equation*}
			\limsup_{\substack{x \to \bar{x}\\z \toattentive{g} c(\bar{x})}} 
			\mathcal{M}(x,z) \subset \mathcal{M}(\bar{x},c(\bar{x})) .
		\end{equation*}
	\end{definition}
\end{mybox}
Sufficient conditions for the validity of the more general qualification condition from \cref{def:AMregularity} can be captured following the arguments of \cite{kruger2021optimality}.
A direct consequence of \cref{prop:minimizers_AM_stationary} is the next result, analogous in spirit to \cite[Proposition 6.9]{kruger2021optimality} and \cite[Corollary 2.7]{demarchi2023constrained}.
\begin{mybox}
	\begin{corollary}\label{cor:AMregular_points}
		Let $x^\ast \in \XX$ be an AKK-regular AKKT-stationary point for \eqref{eq:P}.
		Then, $x^\ast$ is a KKT-stationary point for \eqref{eq:P}.
		In particular, each AKKT-regular local minimizer for \eqref{eq:P} is KKT-stationary.
	\end{corollary}
\end{mybox}

\section{Augmented Lagrangian Method}\label{sec:ALM}
Algorithms exhibiting a nested structure naturally arise in the AL framework, which builds upon a sequence of subproblems.
We develop our implicit approach and derive the generalized AL function in \cref{sec:approach}.
Then, in \cref{sec:upsilonstationarity} we discuss stationarity notions relevant for the AL subproblem, in order to precisely qualify a subproblem ``solution'', and thus the output expected from a subsolver.
The following \cref{sec:ALM:algorithm} contains a detailed statement of our algorithmic framework, whose convergence analysis is presented in \cref{sec:ALM:convergence}.
Suitable termination criteria are discussed in \cref{sec:ALM:termination}.
The subsequent \cref{sec:inner} is dedicated to the numerical solution of the AL subproblems.

	\subsection{Implicit Augmented Lagrangian}\label{sec:approach}
	Introducing an auxiliary variable $z \in \YY$, \eqref{eq:P} can be equivalently rewritten in the form \eqref{eq:ALMz:P},
which has a structured objective function and some (explicit) equality constraints.
For some penalty parameter $\mu>0$, let us define the AL function $\func{\LLslack_\mu}{\XX \times \YY \times \YY}{\Rinf}$ associated to \eqref{eq:ALMz:P} by
\begin{multline}
	\LLslack_\mu(x,z,y)
	{}\coloneqq{}
	f(x) + g(z) + \innprod{y}{c(x) - z} + \frac{1}{2\mu} \| c(x) - z \|^2 \\
	=
	f(x) + g(z) + \frac{1}{2\mu} \| c(x) + \mu y - z \|^2 - \frac{\mu}{2} \|y\|^2 .
	\label{eq:Lz}
\end{multline}
Within AL schemes \cite{bertsekas1996constrained,conn1991globally,birgin2014practical}, it is required to minimize (at least, up to approximate stationarity) the AL function.
Entailing the sum of a smooth and a prox-friendly function, the minimization of $\LLslack_\mu$ can be tackled by exploiting its structure and the computational oracles via, e.g., proximal-gradient methods \cite{combettes2011proximal,parikh2014proximal,themelis2018forward}.
Owing to the convenient structure of $\LLslack_\mu$, the method developed in \cite{demarchi2023constrained} relies on this ``explicit'' approach.
Here, instead of performing a joint minimization of $\LLslack_\mu(\cdot,\cdot,y)$ over both $(x,z)$, we exploit the structure arising from the original problem \eqref{eq:P}, formally eliminating the slack variable $z$, on the vein of the proximal AL approach \cite{dhingra2019proximal,demarchi2021dissertation,hermans2022qpalm,demarchi2022qpdo}.
Given some $\mu > 0$ and $y \in \YY$, we define the oracle $\Sproj_\mu(\cdot,y)$ as the mapping associated to the pointwise minimization of $\LLslack_\mu(\cdot,\cdot,y)$ with respect to $z$.
Owing to the definitions in \eqref{eq:proxmapping} and \eqref{eq:Lz}, the set-valued mapping $\ffunc{\Sproj_\mu(\cdot,y)}{\XX}{\YY}$ is given by
\begin{multline}
	\Sproj_\mu(x,y)
	{}\coloneqq{}
	\argmin_{z\in\YY} \LLslack_\mu(x,z,y)
	{}={}
	\argmin_{z\in\YY} \left\{ g(z) + \frac{1}{2\mu} \| c(x) + \mu y - z \|^2 \right\} \\
	{}={}
	\prox_{\mu g} \left( c(x) + \mu y \right) .
	\label{eq:z}
\end{multline}
Injecting back into $\LLslack_\mu(x,\cdot,y)$ any arbitrary element of $\Sproj_\mu(x,y)$, namely evaluating the AL on the set corresponding to the minimization over $z$, we obtain the (single-valued) AL function $\func{\LL_\mu}{\XX \times \YY}{\R}$ associated to \eqref{eq:P}:
\begin{multline}
	\LL_\mu(x,y)
	{}\coloneqq{}
	\inf_{z\in\YY} \LLslack_\mu(x,z,y)
	{}={}
	\LLslack_\mu(x,z,y)\Big|_{z\in\Sproj_\mu(x,y)} \\
	{}={}
	f(x) + g^\mu( c(x) + \mu y ) - \frac{\mu}{2} \|y\|^2 .
	\label{eq:L}
\end{multline}
We highlight that the Moreau envelope $\func{g^\mu}{\YY}{\R}$ of $g$ is real-valued and strictly continuous \cite{moreau1965proximite}, \cite[Example 10.32]{rockafellar1998variational}, but not continuously differentiable in general, as the proximal mapping of $g$ is possibly set-valued.
Regarded as the AL for \eqref{eq:P}, $\LL_\mu$ was called \emph{generalized} AL function in \cite{rockafellar2022convergence} (for $g$ convex).
Compared to $\LLslack_\mu$ with the \emph{explicit} slack variable $z$ for the lifted ``constrained'' formulation \eqref{eq:ALMz:P}, we may also refer to $\LL_\mu$ as to the \emph{implicit} AL for \eqref{eq:P}.

The algorithm we are going to present in \cref{sec:ALM:algorithm} fits in the AL framework \cite{bertsekas1996constrained,birgin2014practical,rockafellar2022convergence} in that it generates sequences of iterates $x^k$ and Lagrange multipliers $y^k$ obeying the recursion
\begin{equation}
	x^{k+1}
	\leftsquigarrow
	\argmin_{x\in\XX} \LL_{\mu_k}(x,\hat{y}^k) ,
	\qquad
	y^{k+1}
	\leftsquigarrow
	\hat{y}^k + \frac{1}{\mu_k} \partial_y \LL_{\mu_k}(x^{k+1}, \hat{y}^k) ,
	\label{eq:simpleALM}
\end{equation}
with respect to suitable sequences of parameter values $\mu_k > 0$ and multiplier estimates $\hat{y}^k \in \YY$.
Here we intend to give the reader only a prospect about the algorithm before looking into the details; by using the symbol ``$\leftsquigarrow$'' in \eqref{eq:simpleALM} we emphasize that the notation should not be taken too seriously.
Most especially, for the primal update in \eqref{eq:simpleALM} we do not require a \emph{global}, nor an \emph{exact} local minimizer of $\LL_{\mu_k}(\cdot,\hat{y}^k)$, but only an approximate stationary point thereof (in the sense described in \cref{sec:upsilonstationarity}).
Regarding the multiplier, the update rule in \eqref{eq:simpleALM} reduces to the traditional one when \eqref{eq:P} specializes to nonlinear programming \cite{bertsekas1996constrained,birgin2014practical}, but if $g$ is nonconvex it may be necessary to select specific elements from the relevant set-valued subdifferentials.
This potential choice is also very much related to the stationarity concept introduced in \cref{sec:upsilonstationarity}.

At the $k$-th iteration of the AL scheme, a subproblem in the form
\begin{equation}
	\label{eq:minimizeL}
	\minimize_{x\in\XX} \ \LL_{\mu_k}(x,\hat{y}^k)
	,
\end{equation}
with some given $\mu_k > 0$ and $\hat{y}^k \in \YY$, is to be approximately solved.
Being $g$ potentially nonconvex, the AL may lack regularity.
The following \cref{sec:upsilonstationarity} introduces a stationary concept \emph{weaker} than M-stationarity for \eqref{eq:minimizeL} but sharp enough to guarantee that the iterates generated by the AL scheme detailed in \cref{sec:ALM:algorithm} exhibit opportune asymptotic behavior and properties, as demonstrated later in \cref{sec:ALM:convergence}.

	\subsection{Stationarity for Parametric Optimization}\label{sec:upsilonstationarity}
	In view of the marginalization step to derive $\LL_\mu$ in \eqref{eq:L}, we now study the AL subproblem \eqref{eq:minimizeL} with the lenses of parametric optimization \cite[Sections 1.F, 10.C]{rockafellar1998variational}.
Let us consider an objective $\func{\subLL}{\XX}{\R}$ and an oracle $\ffunc{\FBO}{\XX}{\YY}$ given by
\begin{equation}
	\subLL(x) \coloneqq \inf_{z\in\YY} \subLLslack(x,z)
	\qquad\text{and}\qquad
	\FBO(x) \coloneqq \argmin_{z\in\YY} \subLLslack(x,z)
	\label{eq:fbo}
\end{equation}
for a proper, lsc function $\func{\subLLslack}{\XX \times \YY}{\Rinf}$.
We are interested in minimizing $\subLL$, but first in devising suitable stationarity concepts to qualify a ``solution''.
This setting is apparently mirroring subproblem \eqref{eq:minimizeL}, with the explicit and implicit AL functions $\LLslack$ and $\LL$ playing the role of $\subLLslack$ and $\subLL$, respectively.
Moreover, it highlights our intention to keep the solution process implicit, effectively hiding the auxiliary variable $z$, by utilizing the oracle $\FBO$.

Approximate M-stationarity of a point $\bar{x} \in \dom\subLL$ for some given $\varepsilon \geq 0$ amounts to requiring that $\dist(0,\partial \subLL(\bar{x})) \leq \varepsilon$.
However, because of the parametric nature of $\subLL$, the subdifferential $\partial \subLL$ is not a simple object in general.
Let us consider some results on subdifferentiation in parametric optimization.
Recalling from \cref{sec:preliminaries} that the notion of uniform level-boundedness corresponds to a parametric extension of level-boundedness \cite[Definition 1.16]{rockafellar1998variational}, we suppose that $\subLLslack(x,z)$ is level-bounded in $z$ locally uniformly in $x$.
Then, from \cite[Theorem 10.13]{rockafellar1998variational} we have for every $\bar{x}\in\dom\subLL$ the inclusion
\begin{equation}
	\partial \subLL(\bar{x}) {}\subseteq{} \Upsilon(\bar{x}) ,
	\qquad\text{where}\qquad
	\Upsilon(\bar{x}) {}\coloneqq{} \bigcup_{\bar{z} \in \FBO(\bar{x})} \{ v \,\vert\, (v,0) \in \partial \subLLslack(\bar{x},\bar{z}) \} ,
	\label{eq:Upsilon}
\end{equation}
which implies the inequality
\begin{equation}
	\dist(0,\partial \subLL(\bar{x}))
	\geq
	\dist(0, \Upsilon(\bar{x})) .
	\label{eq:subdiffParamOptim}
\end{equation}
This offers a lower bound on the M-stationarity measure, and it does not provide in general an upper bound suitable for a termination condition.
Regardless, in the parametric optimization setting of \eqref{eq:fbo}, we may resort to the following stationarity concept, here denominated \emph{$\Upsilon$-stationarity}.
\begin{mybox}
	\begin{definition}[$\Upsilon$-stationarity]\label{def:UpsilonStationarity}
		Let $\varepsilon \geq 0$ be fixed and consider a proper, lsc function $\func{\subLLslack}{\XX \times \YY}{\Rinf}$ with $\subLLslack(x,z)$ level-bounded in $z$ locally uniformly in $x$.
		Define $\func{\subLL}{\XX}{\Rinf}$ and $\ffunc{\Upsilon}{\XX}{\XX}$ as in \eqref{eq:fbo} and \eqref{eq:Upsilon}, respectively.
		Then, relatively to $\subLL$, a point $x^\ast \in \dom\subLL$ is called \emph{$\varepsilon$-$\Upsilon$-stationary} if $\dist(0, \Upsilon(x^\ast)) \leq \varepsilon$.
		In the case $\varepsilon = 0$, such point $x^\ast$ is said \emph{$\Upsilon$-stationary}.
	\end{definition}
\end{mybox}
Notice that the inclusion $x^\ast \in \dom p$ is implicitly required to have $\Upsilon(x^\ast)$ nonempty.
Furthermore, in the exact case, namely $\varepsilon = 0$, the requirement for $\Upsilon$-stationarity reduces to $0 \in \Upsilon(x^\ast)$, by closedness of $\Upsilon$ \cite[Theorem 10.13]{rockafellar1998variational}.

From the definition of $\Upsilon$ in \eqref{eq:Upsilon}, it appears that an $\varepsilon$-$\Upsilon$-stationary point $x^\ast\in\XX$ can be associated to a (possibly non-unique) \emph{certificate} $z^\ast \in \FBO(x^\ast)$ and a \emph{pair} $(x^\ast, z^\ast)$ that satisfy
\begin{multline}
	\dist(0, \Upsilon(x^\ast))
	{}={}
	\min_{v \in \Upsilon(x^\ast)} \| v \|
	{}={}
	\min_{\substack{v\in\XX\\ z \in \FBO(x^\ast)}} \left\{ \| v \| \,\vert\, (v,0) \in \partial \subLLslack(x^\ast,z) \right\} \\
	{}\leq{}
	\min_{v\in\XX} \left\{ \| v \| \,\vert\, (v,0) \in \partial \subLLslack(x^\ast,z^\ast) \right\}
	\leq
	\varepsilon
	.
	\label{eq:UpsilonCertificate}
\end{multline}
If any such pair $(x^\ast, z^\ast)$ is found, the upper bound in \eqref{eq:UpsilonCertificate} certificates the $\varepsilon$-$\Upsilon$-stationarity of $x^\ast$.

The following result establishes that (approximate) M-stationarity provides a stronger qualification than (approximate) $\Upsilon$-stationarity, meaning that the former is a sufficient condition, not necessary in general, for the latter.
The sought implication readily follows from the upper bound in \eqref{eq:subdiffParamOptim} on the $\Upsilon$-stationarity measure.
\begin{mybox}
	\begin{proposition}\label{lem:upsilon:weaker}
		Let $\varepsilon \geq 0$ be fixed.
		Then, relative to the minimization of $\subLL$, any $\varepsilon$-M-stationary point is also $\varepsilon$-$\Upsilon$-stationary.
	\end{proposition}
\end{mybox}
The converse implication does not hold in general, but it does under some (possibly local) convexity assumptions; see, e.g., \cite[Theorem 2.26]{rockafellar1998variational}.
Similar considerations are valid with regard to the following result, which shows that (approximate) $\Upsilon$-stationarity relative to $\subLL$ is a stronger qualification than (approximate) M-stationarity relative to $\subLLslack$.
\begin{mybox}
	\begin{proposition}\label{lem:upsilon:stronger}
		Let $\varepsilon \geq 0$ be fixed and $(x^\ast, z^\ast)$ be any $\varepsilon$-$\Upsilon$-stationary pair associated to $\subLL$.
		Then, relative to the minimization of $\subLLslack$, $(x^\ast, z^\ast)$ is $\varepsilon$-M-stationary.
	\end{proposition}
\end{mybox}
\begin{proof}
	By $\varepsilon$-$\Upsilon$-stationarity relative to $\subLL$, the pair $(x^\ast,z^\ast)$ necessarily satisfies $x^\ast \in \dom\subLL$ and $z^\ast \in \FBO(x^\ast)$.
	Then, exploiting the upper bound in \eqref{eq:UpsilonCertificate}, we can argue that
	\begin{multline*}
		\dist(0,\partial\subLLslack(x^\ast,z^\ast))
		{}={}
		\min_{\substack{v\in\XX\\ u \in \YY}} \left\{ \| (v,u) \| \,\vert\, (v,u) \in \partial \subLLslack(x^\ast,z^\ast) \right\} \\
		{}\leq{}
		\min_{v\in\XX} \left\{ \| v \| \,\vert\, (v,0) \in \partial \subLLslack(x^\ast,z^\ast) \right\}
		{}\leq{}
		\varepsilon ,
	\end{multline*}
	where the first inequality is due to the restriction $u \coloneqq 0$.
	Showing that $(x^\ast,z^\ast)$ satisfies $\dist(0,\partial\subLLslack(x^\ast,z^\ast)) \leq \varepsilon$, the assertion is proved.
	\qedhere
\end{proof}
Notice, however, that in general the statement of \cref{lem:upsilon:stronger} is not valid for all $z^\ast \in \FBO(x^\ast)$, even if $x^\ast$ is $\varepsilon$-$\Upsilon$-stationary for $\subLL$.
In fact, some pairs $(x^\ast,z^\ast)$ may fail to satisfy \eqref{eq:UpsilonCertificate}, and those are not necessarily $\varepsilon$-M-stationary for $\subLLslack$.

At this point it appears clear that, in relation to \eqref{eq:Lz} and \eqref{eq:L}, the oracle $\FBO$ refers to $\Sproj_{\mu}(\cdot,\hat{y})$, for some given $\mu>0$ and $\hat{y}\in\YY$.
Moreover, we shall focus our attention on the special case of $\subLLslack$ having the structure
\begin{subequations}
	\begin{align}
		\subLLslack(x,z) {}\coloneqq{}& \subLLsmooth(x,z) + g(z) ,
		\label{eq:inner:costslack} \\
		\subLLsmooth(x,z) {}\coloneqq{}& f(x) + \frac{1}{2\mu} \| c(x) + \mu \hat{y} - z \|^2
		\label{eq:PGoracle:psi} ,
	\end{align}
	\label{eq:inner:cost}
\end{subequations}
where $\func{\subLLsmooth}{\XX \times \YY}{\R}$ is continuously differentiable.
	
	\subsection{Algorithm}\label{sec:ALM:algorithm}
	This section details an AL method for the solution of generalized programs of the form \eqref{eq:P}.
In this work we focus on a safeguarded AL scheme inspired by \cite[Algorithm 4.1]{birgin2014practical} and investigate its convergence properties.
Compared to the classical AL or multiplier penalty approach for the solution of nonlinear programs \cite{conn1991globally,bertsekas1996constrained}, this variant uses a safeguarded update rule for the Lagrange multipliers and has stronger global convergence properties \cite{kanzow2018augmented}.
Although we restrict our analysis to this specific algorithm, analogous results can be obtained for others with minor changes, e.g., \cite[Algorithm 1]{evens2021neural}.

The overall method is stated in \cref{alg:ALM}.
Although a starting point is not explicitly required, the subproblems at \cref{step:ALM:subproblem} should be solved warm-starting from the current primal estimate $x^{k-1}$, thus exploiting initial guesses.
The AL subproblem appears at \cref{step:ALM:subproblem}, involving the implicit AL function with current penalty parameter $\mu_k$ and estimate $\hat{y}^k$.
Crucially, \cref{step:ALM:subproblem} requires to compute an approximate $\Upsilon$-stationary \emph{pair} $(x^k,z^k)$, not merely a point $x^k$, for $\LL_{\mu_k}(\cdot,\hat{y}^k)$.
As it may be employed to verify approximate stationarity, it seems likely that a certificate $z^k$ for $x^k$ can be obtained with no additional effort --- this is the case for the subsolver developed in \cref{sec:inner}.

\begin{algorithm2e}[htb]
	\DontPrintSemicolon
	\KwData{$\epsilon > 0$}
	\KwResult{$\epsilon$-KKT-stationary point $x^\ast$}%
	Select $\mu_0 > 0$, $\theta, \kappa \in (0,1)$, $\Ybounded\subset\YY$ nonempty bounded\;
	\For{$k = 0,1,2\ldots$}{
		Select $\hat{y}^k \in \Ybounded$ and $\varepsilon_k \geq 0$\label{step:ALM:ysafe}\; 
		Compute an $\varepsilon_k$-$\Upsilon$-stationary pair $(x^k, z^k)$ for $\LL_{\mu_k}(\cdot,\hat{y}^k)$\label{step:ALM:subproblem}\label{step:ALM:x}\label{step:ALM:z}\;
		Set 
		\label{step:ALM:y}%
		$y^k \gets \hat{y}^k + \mu_k^{-1} [c(x^k) - z^k]$
		and $V^k \gets \| c(x^k) - z^k \|$\;
		\lIf{$\max\{ \varepsilon_k, V^k \} \leq \epsilon$\label{step:ALM:termination}}{%
			\Return $x^\ast \gets x^k$%
		}%
		\lIf{$k=0$ \KwOr $ V^k \leq \theta V^{k-1}$\label{step:ALM:penalty}}{
			$\mu_{k+1} \gets \mu_k$, \textbf{else} $\mu_{k+1} \gets \kappa \mu_k$
		}
	}
	\caption{Safeguarded implicit AL method for \eqref{eq:P}}
	\label{alg:ALM}
\end{algorithm2e}%

The safeguarded Lagrange multiplier estimate $\hat{y}^k$ is drawn from a bounded set $\Ybounded \subset \YY$ at \cref{step:ALM:ysafe}.
Although not necessary, the choice of $\hat{y}^k$ should also depend on the current estimate $y^{k-1}$.
Moreover, the choice of $\Ybounded$ can take advantage of \emph{a priori} knowledge of $\partial g$ and its structure, in order to generate better Lagrange multiplier estimates.
In practice, it is advisable to choose the safeguarded estimate $\hat{y}^k$ as the projection of the Lagrange multiplier $y^{k-1}$ onto $\Ybounded$, thus effectively adopting the classical approach as long as $y^{k-1}$ remains within $\Ybounded$; we refer to \cite[Section 4.1]{birgin2014practical} and \cite[Section 3]{kanzow2018augmented} for a more articulated discussion.

Inspired by the classical first-order Lagrange multiplier estimate, the update rule at \cref{step:ALM:y} is designed around \eqref{eq:AMstationary:x} and the identity
\begin{equation*}
	\nabla_x \LLslack_\mu(x,z,\hat{y})
	=
	\nabla f(x) + \frac{1}{\mu} \jac c(x)^\top [c(x) + \mu \hat{y} - z] .
\end{equation*}
Given any $(x^k,z^k)$ obtained at \cref{step:ALM:subproblem}, the specific choice of $y^k$ leads by construction to
\begin{equation*}
	\| \nabla f(x^k) + \jac c(x^k)^\top y^k \|
	\leq
	\varepsilon_k
	,
\end{equation*}
effectively connecting inner and outer loops, and providing an upper bound on the dual infeasibility at $(x^k,y^k)$.
\Cref{step:ALM:termination} encapsulates suitable termination conditions for returning an approximate KKT-stationary point; these are discussed in \cref{sec:ALM:termination}.
The monotonicity test at \cref{step:ALM:penalty} is adopted to monitor primal infeasibility along the iterates and update the penalty parameter accordingly.
Aimed at driving $V^k \coloneqq \|c(x^k) - z^k\|$ to zero, the penalty parameter $\mu_k$ is reduced at \cref{step:ALM:penalty} in case of insufficient decrease.
	
	\subsection{Convergence Analysis}\label{sec:ALM:convergence}
	Throughout our convergence analysis, we assume that \cref{alg:ALM} is well-defined, thus requiring that each subproblem at \cref{step:ALM:subproblem} admits an approximate stationary point.
Moreover, the following statements assume the existence of some accumulation point $x^\ast$ for a sequence $\{ x^k \}_{k\in\N}$ generated by \cref{alg:ALM}.
In general, coercivity or (level) boundedness arguments should be adopted to verify these preconditions.

Due to their practical importance, we focus on \emph{local}, or affordable, solvers, which return stationary points as candidate local minimizers, for the subproblems at \cref{step:ALM:subproblem}.
Instead, we do not investigate the case where the subproblems are solved to \emph{global} optimality, whose analysis would pattern the classical results in \cite[Chapter 5]{birgin2014practical} and \cite{kanzow2018augmented}.

As it is often the case for penalty-type methods in the nonconvex setting, \cref{alg:ALM} may generate accumulation points that are infeasible for \eqref{eq:P}.
Patterning standard arguments, the following result gives conditions that guarantee feasibility of limit points; cf. \cite[Problem 4.12]{birgin2012augmented}, \cite[Proposition 3.2]{demarchi2023constrained}.
\begin{mybox}
	\begin{proposition}\label{lem:ALM:bounded}
		Let \cref{ass:P} hold and consider a sequence $\{ x^k \}_{k\in\N}$ of iterates generated by \cref{alg:ALM} with $\epsilon = 0$.
		Then, each accumulation point $x^\ast$ of $\{ x^k \}_{k\in\N}$ is feasible for \eqref{eq:P} if (at least) one of the following conditions holds:
		\begin{enumerate}[label=(\roman*)]
			\item\label{lem:ALM:bounded:attentive}%
			$x^k \toattentive{\varphi}_K x^\ast$ for some subsequence $\{x^k\}_{k\in K}$;
			\item\label{lem:ALM:bounded:bound}%
			there exists some $B \in \R$ such that $\LL_{\mu_k}(x^k, \hat{y}^k) \leq B$ for all $k\in\N$;
			\item\label{lem:ALM:bounded:mu}%
			$\{ \mu_k \}_{k\in\N}$ is bounded away from zero and $\dom g$ is closed;
			\item\label{lem:ALM:bounded:muLagr}%
			$\{ \mu_k \}_{k\in\N}$ is bounded away from zero and $\LL_{\mu_k}(x^k,\hat{y}^k) \leq \LL_{\mu_k}(x^{k-1},\hat{y}^k)$ for all $k\in\N$.
		\end{enumerate}
	\end{proposition}
\end{mybox}
\begin{proof}\label{lem:ALM:bounded:proof}
	Let $x^\ast \in \XX$ be an arbitrary accumulation point of $\{ x^k \}$ and 
	$\{x^k\}_{k\in K}$ a subsequence such that $x^k \to_K x^\ast$.
	We shall argue that $x^\ast$ is feasible for \eqref{eq:P}, that is, $c(x^\ast) \in \dom g$.
	Before considering each assertion, we point out that, if $\{ \mu_k \}$ is bounded away from zero, the conditions at \cref{step:ALM:penalty} of \cref{alg:ALM} imply that $V^k \coloneqq \|c(x^k) - z^k\| \to 0$ for $k \to \infty$.
	Thus, since $x^k \to_K x^\ast$, by continuity arguments we have also $z^k \to_K c(x^\ast)$.
	Finally, let us conveniently denote by $\Delta_Y \in \R$ the finite size of the bounded set $Y \subset \YY$ initialized in \cref{alg:ALM}, such that $\| \hat{y} \| \leq \Delta_Y$ for all $\hat{y} \in Y$.
	\begin{enumerate}
		\item[\ref{lem:ALM:bounded:attentive}]
		The $\varphi$-attentive convergence of $\{x^k\}_{k\in K}$ to $x^\ast$ requires that $x^\ast \in \dom \varphi$.
		Hence, by $\varphi \coloneqq g \circ c$ and continuity of $c$, it must be $c(x^\ast) \in \dom g$.
		\item[\ref{lem:ALM:bounded:bound}]
		In view of the preliminary observation regarding $\{ \mu_k \}$ bounded away from zero, it remains to consider the case $\mu_k \to 0$ for showing that $z^k \to_K c(x^\ast)$.
		Using \eqref{eq:L} and $z^k \in \Sproj_{\mu_k}(x^k,\hat{y}^k)$, the upper bound yields
		\begin{equation*}
			B
			\geq
			\LL_{\mu_k}(x^k,\hat{y}^k)
			=
			f(x^k) + g(z^k) + \frac{1}{2 \mu_k} \| c(x^k) + \mu_k \hat{y}^k - z^k \|^2 - \frac{\mu_k}{2} \| \hat{y}^k \|^2
		\end{equation*}
		for all $k\in\N$.
		Rearranging terms and considering $\mu_k \to 0$, we get
		\begin{equation*}
			\| c(x^k) + \mu_k \hat{y}^k - z^k \|^2
			\leq 
			2 \mu_k \left[ B + \frac{\mu_k}{2} \| \hat{y}^k \|^2 - f(x^k) - g(z^k) \right]
			\to_K
			0 ,
		\end{equation*}
		where the limit follows from boundedness of $\{\hat{y}^k\}_{k\in\N}$ and $\{z^k\}_{k\in K}$, the latter remaining bounded by $x^k \to_K x^\ast$ and local boundedness of the proximal mapping.
		Using once again that $\mu_k \hat{y}^k \to 0$, the limit implies $c(x^k) - z^k \to_K 0$.
		Thus, since $x^k \to_K x^\ast$, by continuity arguments it is $z^k \to_K c(x^\ast)$.
		Now, owing to $x^k\to_Kx^\ast$, boundedness of $\{\hat{y}^k\}_{k\in\N}$, and $\mu_k \leq \mu_0$, from the upper bound condition we also obtain that $\{g(z^k)\}_{k\in K}$ remains bounded above:
		\begin{equation*}
			\limsup_{k\to_K\infty} g(z^k)
			\leq
			B  + \limsup_{k\to_K\infty} \frac{\mu_k}{2} \| \hat{y}^k \|^2 - f(x^k)
			\leq
			B + \frac{\mu_0}{2} \Delta_Y^2 - f(x^\ast)
			<
			\infty .
		\end{equation*}
		Combined with $z^k \to_K c(x^\ast)$ and lower semicontinuity of $g$, this implies that $c(x^\ast) \in \dom g$.
		\item[\ref{lem:ALM:bounded:mu}]
		As $\{\mu_k\}_{k\in\N}$ remains bounded away from zero, we have that $z^k \to_K c(x^\ast)$.
		Then, since $z^k \in \Sproj_{\mu_k}(x^k,\hat{y}^k) \subseteq \dom g$ for all $k\in\N$ by \cref{step:ALM:z}, it follows from $\dom g$ being closed that $c(x^\ast) \in \dom g$.
		\item[\ref{lem:ALM:bounded:muLagr}]
		As $\{\mu_k\}_{k\in\N}$ remains bounded away from zero, we have that $z^k \to_K c(x^\ast)$ and, owing to the conditions at \cref{step:ALM:penalty}, $\{\mu_k\}_{k\in\N}$ is asymptotically constant, namely for some finite index $k_\ast\in\N$ it is $\mu_k = \mu_\ast > 0$ for all $k \geq k_\ast$.
		In turn, the primal infeasibility measure $V^k$ decreases linearly, as it must be $V^k \leq \theta V^{k-1}$ for all $k \geq k_\ast$.
		Now, rearranging \eqref{eq:L} and exploiting the explicit minimization over $z$, we obtain the upper bound
		\begin{align*}
			\LL_{\mu_\ast}(x^{k-1}, \hat{y}^k)
			{}={}&
			f(x^{k-1}) + \min_z \left\{ g(z) + \innprod{\hat{y}^k}{c(x^{k-1}) - z} + \frac{1}{2 \mu_\ast} \| c(x^{k-1}) - z \|^2 \right\} \\
			{}\leq{}&
			f(x^{k-1}) + g(z^{k-1}) + \innprod{\hat{y}^k}{c(x^{k-1}) - z^{k-1}} + \frac{1}{2 \mu_\ast} \| c(x^{k-1}) - z^{k-1} \|^2 \\
			{}={}&
			\LL_{\mu_\ast}(x^{k-1}, \hat{y}^{k-1}) + \innprod{\hat{y}^k - \hat{y}^{k-1}}{c(x^{k-1}) - z^{k-1}}
		\end{align*}
		for all $k > k_\ast$.
		Invoking the Cauchy-Schwarz inequality, the descent condition gives
		\begin{align*}
			\LL_{\mu_\ast}(x^k, \hat{y}^k)
			{}\leq{}
			\LL_{\mu_\ast}(x^{k-1}, \hat{y}^k)
			{}\leq{}&
			\LL_{\mu_\ast}(x^{k-1}, \hat{y}^{k-1}) + 2 \Delta_Y V^{k-1}
		\end{align*}
		for all $k > k_\ast$, where we used that $\| \hat{y}^k - \hat{y}^{k-1} \| \leq 2 \Delta_Y$.
		Then, telescoping this argument, the linear decrease of $V^k$ yields
		\begin{equation*}
			\LL_{\mu_\ast}(x^{k+\ell}, \hat{y}^{k+\ell}) - \LL_{\mu_\ast}(x^{k-1}, \hat{y}^{k-1})
			{}\leq{}
			2 \Delta_Y \sum_{m=0}^{\ell} V^{k+m-1}
			{}\leq{}
			2 \Delta_Y V^{k-1} \sum_{m=0}^{\ell} \theta^\ell
		\end{equation*}
		for all $k > k_\ast$ and $\ell \in \N$.
		Hence, since $\theta \in (0,1)$, for any arbitrary but fixed $k > k_\ast$ taking the limit results in
		\begin{equation*}
			\limsup_{\ell \to \infty}
			\LL_{\mu_\ast}(x^{k+\ell}, \hat{y}^{k+\ell})
			\leq
			\LL_{\mu_\ast}(x^{k-1}, \hat{y}^{k-1}) + \frac{2 \Delta_Y}{1 - \theta} V^{k-1}
			<
			\infty
			,
		\end{equation*}
		showing that, up to discarding early iterates, there is an upper bound to $\LL_{\mu_k}(x^k, \hat{y}^k)$.
		Then, feasibility of $x^\ast$ follows as in \ref{lem:ALM:bounded:bound}.
		\qedhere
	\end{enumerate}
\end{proof}

The following result concerns the Lagrange multiplier estimates.
Considering the particular case of $g$ being the indicator of some set $D$ with nonempty interior, \cref{lem:ALM:multipliers} shows that the Lagrange multipliers vanish for constraints that are (strictly) inactive in the limit.
This indicates that the update rule at \cref{step:ALM:y} of \cref{alg:ALM} is plausible with regard to their expected behavior as shifts \cite[Section 4]{birgin2014practical}.
\begin{mybox}
	\begin{proposition}\label{lem:ALM:multipliers}
		Suppose $g \equiv \indicator_D$ for some nonempty and closed set $D \subseteq \YY$.
		Let \cref{ass:P} hold and consider a sequence $\{ x^k \}_{k\in\N}$ of iterates 
		generated by \cref{alg:ALM} with $\epsilon = 0$.
		Let $x^\ast$ be an accumulation point of $\{ x^k \}_{k\in\N}$ and 
		$\{x^k\}_{k\in K}$ a subsequence such that $x^k \to_K x^\ast$.
		If $c(x^\ast) \in \interior D$, then $y^k = 0$ for all $k \in K$ large enough.
	\end{proposition}
\end{mybox}
\begin{proof}
	Let $c(x^\ast) \in \interior D \ne \emptyset$ and notice that the proximal mapping of $\indicator_D$ coincides with the projection onto $D$, namely $\prox_{\gamma \indicator_D} = \proj_D$ for all $\gamma > 0$.
	Consider the two cases.
	\begin{itemize}
		\item If $\mu_k \to 0$, by boundedness of $\{ \hat{y}^k \}_{k\in\N}$ and continuity of $c$, there exists $k_1 \in K$ such that $c(x^k) + \mu_k \hat{y}^k \in \interior D$ for all $k \in K$, $k \geq k_1$.
		\item If $\{ \mu_k \}_{k\in\N}$ is bounded away from zero, it follows from \cref{step:ALM:penalty} that $c(x^k) - z^k \to_K 0$.
		Then, by continuity of $c$ and $x^k \to_K x^\ast$, $z^k \to_K c(x^\ast) \in \interior D$ follows.
		Hence, there exists $k_2 \in K$ such that $z^k \in \interior D$ for all $k \in K$, $k \geq k_2$.
	\end{itemize}
	In both cases, since $z^k \in \Sproj_{\mu_k}(x^k, \hat{y}^k) = \proj_D( c(x^k) + \mu_k \hat{y}^k )$ and $z^k \in \interior D$ for all $k \in K$ large enough, it must be $z^k = c(x^k) + \mu_k \hat{y}^k$ and, consequently, $y^k = 0$ by \cref{step:ALM:y}.
	\qedhere
\end{proof}
Notice that the assertion can be easily refined by exploiting the separable structure of $D$, if any.
For instance, if $D$ is a hyperbox, every step in the proof can be applied componentwise, recovering \cite[Theorem 4.1]{birgin2014practical}.

The following convergence result provides fundamental theoretical support to \cref{alg:ALM}.
It shows that, under subsequential attentive convergence, accumulation points are AKKT-stationary for \eqref{eq:P}.
For any accumulation point $x^\ast$ and for any subsequence $\{x^k\}_{k\in K}$ such that $x^k \to_K x^\ast$, \cref{lem:ALM:bounded}\ref{lem:ALM:bounded:attentive} assures the feasibility of $x^\ast$ under the assumption that $x^k \toattentive{\varphi}_K x^\ast$.
However, the $\varphi$-attentive convergence also entails $g(c^k) \to_K g(c(x^\ast))$ as $x^k \to_K x^\ast$, which is trivially satisfied, e.g., if ($x^\ast$ is feasible and) $g$ is continuous relative to $\dom g$, in the sense that, for all $x\in\dom g$ and $\{x^\ell\} \subseteq \dom g$, whenever $x^\ell \to x$ it holds that $g(x^\ell) \to g(x)$.
We refer to \cite[Example 3.4]{demarchi2023constrained} for an illustration of the importance of $\varphi$-attentive convergence.
\begin{mybox}
	\begin{theorem}\label{thm:AMstationary}
		Let \cref{ass:P} hold and consider a sequence $\{ x^k \}_{k\in\N}$ of iterates generated by \cref{alg:ALM} with $\epsilon = 0$ and $\varepsilon_k \to 0$.
		Let $x^\ast\in\XX$ be an accumulation point of $\{ x^k \}_{k\in\N}$ and $\{x^k\}_{k\in K}$ a subsequence such that $x^k \toattentive{\varphi}_K x^\ast$.
		Then, $x^\ast$ is an AKKT-stationary point for \eqref{eq:P}.
	\end{theorem}
\end{mybox}
\begin{proof}
	We claim that the subsequences $\{x^k\}_{k\in K}$, $\{y^k\}_{k\in K}$ and $\{z^k\}_{k\in K}$ generated by \cref{alg:ALM} satisfy the properties in \cref{def:AMstationary} and therefore show that $x^\ast$ is an AKKT-stationary point for \eqref{eq:P}.
	First of all, by \cref{lem:ALM:bounded}\ref{lem:ALM:bounded:attentive}, the $\varphi$-attentive convergence of $\{x^k\}_{k\in K}$ implies that $x^\ast$ is feasible for \eqref{eq:P}.
	From \cref{step:ALM:subproblem,step:ALM:y} of \cref{alg:ALM}, we have by $\varepsilon_k$-$\Upsilon$-stationarity of $x^k$ with certificate $z^k$ that
	\begin{equation*}
		\eta^k
		\coloneqq
		\nabla f(x^k) + \jac c(x^k)^\top y^k
		=
		\nabla_x \LLslack_{\mu_k}(x^k,z^k,\hat{y}^k)
		\label{thm:AMstationary:x}
	\end{equation*}
	satisfies $\| \eta^k \| \leq \varepsilon_k$.
	Then, by construction, we have $\eta^k \to_K 0$.
	Further, from \cref{step:ALM:y} we obtain that, for all $k \in \N$,
	\begin{equation*}
		y^k 
		=
		\hat{y}^k + \frac{c(x^k) - z^k}{\mu_k}
		\in
		\partial g( z^k ) ,
		\label{thm:AMstationary:y}
	\end{equation*}
	where the inclusion follows from $z^k \in \Sproj_{\mu_k}(x^k, \hat{y}^k)$ (in fact, from the necessary condition for the optimality of $z^k$).
	It remains to show that $c(x^k) - z^k \to_K 0$.
	To this end, we consider the two cases.
	\begin{itemize}
		\item
		If $\{\mu_k\}_{k\in\N}$ is bounded away from zero, the conditions at \cref{step:ALM:penalty} imply that $V^k \coloneqq \| c(x^k) - z^k \| \to 0$, hence in particular $c(x^k) - z^k \to_K 0$.
		\item
		If instead $\mu_k \to 0$, we exploit continuity of $c$ and boundedness of $\{\hat{y}^k\}_{k\in\N}$.
		These properties result in $c(x^k) + \mu_k \hat{y}^k \to_K c(x^\ast) \in \dom g$, where the inclusion follows from feasibility of $x^\ast$.
		Then, since $z^k \in \Sproj_{\mu_k}(x^k,\hat{y}^k)$ for all $k\in\N$, from \eqref{eq:z} and \cite[Lemma 3.3]{demarchi2022interior} we have that
		\begin{equation*}
			z^k
			\in
			\prox_{\mu_k g}(c(x^k) + \mu_k \hat{y}^k)
			\to_K
			c(x^\ast)
		\end{equation*}
		while $\mu_k \to 0$.
		Thus, the limit $c(x^k) - z^k \to_K 0$ follows.
	\end{itemize}
	Overall, this proves that $x^\ast$ is an AKKT-stationary point for \eqref{eq:P}.
	\qedhere
\end{proof}

Constrained optimization algorithms aim at finding feasible points and minimizing the objective function subject to constraints.
Employing affordable local optimization techniques, one cannot expect to find global minimizers of any infeasibility measure.
Nevertheless, the next result proves that \cref{alg:ALM} with bounded $\{ \varepsilon_k \}_{k\in\N}$ finds at least stationary points of an infeasibility measure.
Notice that this property does not require $\varepsilon_k \to 0$, but only boundedness; cf. \cite[Theorem 6.3]{birgin2014practical}.
\begin{mybox}
	\begin{proposition} \label{lem:feasP}
		Let \cref{ass:P} hold and suppose $\dom g$ is closed.
		Consider a sequence $\{ x^k \}_{k\in\N}$ of iterates generated by \cref{alg:ALM} with $\epsilon = 0$ and $\{ \varepsilon_k \}_{k\in\N}$ bounded.
		Let $x^\ast \in \XX$ be an accumulation point of $\{ x^k \}_{k\in\N}$ and $\{x^k\}_{k\in K}$ a subsequence such that $x^k \to_K x^\ast$.
		Then, $x^\ast$ is an M-stationary point of
		\begin{equation*}
			\minimize_{x\in\XX} \,
			\dist^2( c(x), \dom g ) .
		\end{equation*}
	\end{proposition}
\end{mybox}
\begin{proof}
	By \cref{lem:ALM:bounded}\ref{lem:ALM:bounded:mu} and $\dom g$ being closed, if $\{\mu_k\}_{k\in\N}$ is bounded away from zero, then each accumulation point $x^\ast$ is feasible for \eqref{eq:P}, namely it is a global minimizer of the feasibility problem, and so it is an M-stationary point thereof by strict continuity of the objective function \cite[Proposition 5.3]{mordukhovich2006variational}.
	Hence, it remains to consider the case $\mu_k \to 0$.

	Let us rewrite the feasibility problem with the equivalent lifted representation
	\begin{equation}
		\minimize_{\substack{x\in\XX \\ (z,\alpha)\in\epi g}}
		{}\quad{}
		\frac{1}{2} \| c(x) - z \|^2
		\label{eq:feasPzEpigraph}
	\end{equation}
	obtained by epigraphical reformulation.
	The inclusion $(z,\alpha) \in \epi g$ requires that $g(z) \leq \alpha < \infty$, thus encoding $z \in \dom g$.
	M-stationarity of $x^\ast$ for \eqref{eq:feasPzEpigraph} requires the existence of some $z^\ast \in \YY$ such that
	\begin{subequations}\label{eq:feasP:KKT}
		\begin{align}
			\jac c(x^\ast)^\top [c(x^\ast) - z^\ast] {}={}& 0 , \label{eq:feasP:KKTx}\\
			(c(x^\ast) - z^\ast, 0) {}\in{}& \limnormalcone_{\epi g}(z^\ast, g(z^\ast)) \label{eq:feasP:KKTz}.
		\end{align}
	\end{subequations}
	Owing to \cref{step:ALM:subproblem} of \cref{alg:ALM} it is $\|\eta^k\| \leq \varepsilon_k$, where
	\begin{equation*}
		\eta^k 
		\coloneqq
		\nabla f(x^k) 
		+ 
		\frac{1}{\mu_k} \jac c(x^k)^\top \left[ c(x^k) + \mu_k \hat{y}^k - z^k \right]
		,
	\end{equation*}
	for all $k \in K$.
	Taking the limit of $\mu_k \eta^k$ for $k \to_K \infty$, while using $\mu_k \to 0$ and boundedness of $\{\eta^k\}_{k\in\N}$ and $\{\nabla f(x^k)\}_{k\in K}$, it follows that \eqref{eq:feasP:KKTx} is satisfied.
	Furthermore, by the minimizing property of $z^k {}\in{} \Sproj_{\mu_k}(x^k, \hat{y}^k) = \prox_{\mu_k g}(c(x^k) + \mu_k \hat{y}^k)$, it is $z^k \in \dom g$ and 
	\begin{equation*}
		\frac{1}{\mu_k} \left[ c(x^k) + \mu_k \hat{y}^k - z^k \right]
		\in
		\partial g(z^k)
	\end{equation*}
	for all $k\in\N$.
	By boundedness of $\{c(x^k)+\mu_k\hat y^k\}_{k\in K}$ and local boundedness of the proximal mapping, $\{z^k\}_{k\in K}$ remains bounded too.
	Hence, since $z^k \in \dom g$ for all $k\in\N$ and $\dom g$ is closed, we may assume without loss of generality that $z^k \to_K z^\ast$ for some $z^\ast \in \dom g$.
	Using the relationship stated in \eqref{eq:subdiffLimNormalConeEpigraph}, multiplying by $\mu_k > 0$ and exploiting that $\limnormalcone_{\epi g}(z^k,g(z^k))$ is a cone, we have
	\begin{equation*}
		(c(x^k) + \mu_k \hat{y}^k - z^k, -\mu_k)
		\in 
		\limnormalcone_{\epi g}(z^k, g(z^k))
	\end{equation*}
	for all $k\in\N$.
	Taking the limit for $k\to_K\infty$, the robustness of the limiting normal cone \cite[Proposition 1.20]{mordukhovich2006variational} and the convergence $z^k \to_K z^\ast \in \dom g$ yield \eqref{eq:feasP:KKTz}, concluding the proof.
	\qedhere
\end{proof}

	\subsection{Termination Criteria}\label{sec:ALM:termination}
	\Cref{step:ALM:subproblem} of \cref{alg:ALM} involves the minimization of the implicit AL function, defined in \eqref{eq:L}.
Then, the multiplier update at \cref{step:ALM:y} allows one to draw conclusions with respect to the original problem \eqref{eq:P}, as shown by \cref{thm:AMstationary}.
Moreover, recalling the AKKT-stationarity conditions \eqref{eq:AMstationary}, in view of \cref{step:ALM:ysafe} one may construct a sequence $\{ \varepsilon_k \}_{k\in\N} \subset \R_{++}$ such that $\varepsilon_k \to 0$.
Then, given some user-specified tolerances $\varepsilon^{\text{dual}}, \varepsilon^{\text{prim}} > 0$, it is reasonable to declare successful convergence when the conditions $\varepsilon_k \leq \varepsilon^{\text{dual}}$ and $V^k \coloneqq \| c(x^k) - z^k \| \leq \varepsilon^{\text{prim}}$ are satisfied.
\Cref{thm:AMstationary} demonstrates that these termination criteria (the latter, in particular) are satisfied in finitely many iterations if any subsequence of $\{x^k\}_{k\in K}$ accumulates at a feasible point $x^\ast \in \XX$.
As this might not be the case, a mechanism for (local) infeasibility detection is needed, and usually included in practical implementations.

Given some tolerances, \cref{alg:ALM} can be equipped with relaxed conditions on decrease requirements at \cref{step:ALM:penalty} and optimality at \cref{step:ALM:subproblem} \cite{birgin2012augmented}.
At \cref{step:ALM:ysafe} the inner tolerance $\varepsilon_k$ can stay bounded away from zero, as long as $\varepsilon_k \leq \varepsilon^{\text{dual}}$ for large $k \in \N$.
Similarly, the condition at \cref{step:ALM:penalty} can be relaxed by adding the (inclusive) possibility that $V^k \leq \varepsilon^{\text{prim}}$, thus avoiding unnecessarily strong penalization.

\section{Inner Problem and Solver}\label{sec:inner}
In this section we elaborate upon the minimization of the implicit AL function, defined in \eqref{eq:L}, focusing on the numerical solution of subproblems \eqref{eq:minimizeL} for some fixed $\mu_k > 0$ and $\hat{y}^k \in \YY$.
As remarked in \cref{sec:approach}, this task could be performed by jointly minimizing the explicit AL function $\LLslack_{\mu_k}(\cdot,\cdot,\hat{y}^k)$ in \eqref{eq:Lz}, that is, with respect to both primal and slack variables.
Owing to its composite structure, any suitable methods for composite optimization could be adopted for this task; as pointed out in \cite[Section 3.4]{demarchi2023constrained}, a prominent example are proximal gradient algorithms.

We take a different approach, aimed at minimizing $\LL_{\mu_k}(\cdot,\hat{y}^k)$ directly, in order to solve inner problems without auxiliary variables.
To do so, we juggle the implicit and explicit formulations, linked by the oracle, as in \eqref{eq:fbo}.
In \cref{sec:PGoracle}, an algorithm fo the minimization of $\subLL$ is discussed, that exploits the available oracles and the structure of $\subLLslack$ in \eqref{eq:inner:cost}.
Reminiscent of \cite[Section 1.5.4]{demarchi2021dissertation}, the scheme compensates for the lack of regularity by pivoting around the oracle's optimality, as a kind of reliable \emph{magical step} \cite{conn2000trust,birgin2014practical}.
In fact, a surrogate gradient of $\subLL$ is built based on that of $\subLLsmooth$, the smooth part of $\subLLslack$, evaluated at the oracle's response.
Then, with a bit of abstraction, it turns out that any off-the-shelf algorithm for unconstrained smooth optimization can be directly adopted for the task of minimizing the implicit AL function, in a remarkably transparent way.

Throughout this section, let $\mu_k > 0$ and $\hat{y}^k \in \YY$ be some fixed given parameters.  Denote $\func{\subLL}{\XX}{\R}$ and $\func{\subLLsmooth}{\XX \times \YY}{\R}$ as the implicit AL and its smooth part, respectively.
Also, let $\func{\subLLslack}{\XX \times \YY}{\Rinf}$ be the explicit AL; see \cref{sec:upsilonstationarity}.
Patterning the transition from $\LLslack_\mu(\cdot,\cdot,y)$ in \eqref{eq:Lz} to $\LL_\mu(\cdot,y)$ in \eqref{eq:L}, we suppose that an oracle $\FBO$ is available that performs the explicit (exact, global) minimization of $\subLLslack$ with respect to $z \in \YY$, apparently alluding to \eqref{eq:z}.
In relation to the $k$th iteration our AL method, we write $\FBO = \Sproj_{\mu_k}(\cdot, \hat{y}^k)$.
Notice that, for all $x \in \XX$, it is $\FBO(x) \subseteq \dom g$ and thus $\subLL(x) < \infty$.
We collect some basic properties under \cref{ass:P}.
\begin{mybox}
	\begin{proposition}\label{thm:PGoracle:properties}%
		Consider \eqref{eq:fbo}, \eqref{eq:inner:cost}, and let \cref{ass:P} hold.
		Then,
		\begin{enumerate}[label=(\roman*)]
			\item\label{ass:PGoracle:subLLsmooth}
			$\func{\subLLsmooth}{\XX \times \YY}{\R}$ is continuously differentiable;
			\item\label{ass:PGoracle:g}
			$\func{g}{\YY}{\Rinf}$ is proper, lower semicontinuous, and prox-bounded;
			\item\label{ass:PGoracle:subLLslack}
			$\func{\subLLslack}{\XX \times \YY}{\Rinf}$, $\subLLslack(x,z)$ is level-bounded in $z$ locally uniformly in $x$;
			\item\label{ass:PGoracle:oracle}
			$\ffunc{\FBO}{\XX}{\YY}$ is locally bounded, nonempty- and compact-valued;
			\item\label{ass:PGoracle:subLL}
			$\func{\subLL}{\XX}{\R}$ is strictly continuous (i.e. locally Lipschitz).
		\end{enumerate}
	\end{proposition}
\end{mybox}
\begin{proof}
	The assertions readily follow from \cref{ass:P}\ref{ass:f}--\ref{ass:g} and the definitions in \eqref{eq:fbo} and \eqref{eq:inner:cost}.
	In particular, \ref{ass:PGoracle:subLLslack} is due to prox-boundedness of $g$ and that $\mu_k \leq \gamma_g$, along with $\hat{y}^k$ being fixed and smoothness of $f$ and $c$.
	Assertion \ref{ass:PGoracle:oracle} is a consequence of the smoothness of $c$ by \cref{ass:P}\ref{ass:f} and of the properties of the proximal mapping \cite[Theorem 1.25]{rockafellar1998variational}.
	Assertion \ref{ass:PGoracle:subLL} follows from the identity $\subLL(\cdot) = f(\cdot) + g^{\mu_k}(c(\cdot) + \mu_k \hat{y}^k)$ and the fact that the Moreau envelope is strictly continuous \cite[Example 10.32]{rockafellar1998variational}.
	\qedhere
\end{proof}
Owing to \cref{ass:P} and \cref{thm:PGoracle:properties}\ref{ass:PGoracle:subLL}, the objective function $\subLL$ is proper and lsc.
However, as the existence of a minimizer for $\subLL$ is not guaranteed a priori, we also consider the following assumption.
In general, it is verified under some coercivity or (level) boundedness conditions.
\begin{mybox}
	\begin{assum}\label{ass:PGoracle}\label{ass:PGoracle:wellposed}
		Consider \eqref{eq:fbo}.
		Then, $\inf \subLL \in \R$.
	\end{assum}
\end{mybox}

\begin{remark}\label{rem:abstractConstraints}
	The problem format \eqref{eq:P} could incorporate also an abstract constraint $x \in X$, with the set $X \subseteq \XX$ being nonempty, closed, and \emph{simple}.
	Although not included for the sake of a simpler exposition, the results in this paper could be readily extended to account for a \emph{convex} $X$.
	As noted in \cite{rockafellar2022augmented}, a possibility is to embed $x \in X$ in the template \eqref{eq:P} by replacing $c(x)$ by $\bar{c}(x) \coloneqq (c(x), x)$ and $g(z)$ by $\bar{g}(z, z^\prime) \coloneqq g(z) + \indicator_X(z^\prime)$.
	This approach, however, introduces an extra slack component $z^\prime$, with an extra multiplier component $y^\prime$, while relaxing the simple constraint $x \in X$.
	Following \cite{andreani2008augmented} instead, a sharper strategy is to treat $x \in X$ directly as an unrelaxable constraint for the AL subproblem, and not incorporating it into the implicit AL function.
	Exploiting the projection onto $X$, the inner solver could generate iterates that satisfy the abstract constraint by design.
	With potentially nonconvex $X$, though, it appears nontrivial to formulate a suitable termination criterion for attesting approximate stationarity.

	The possibility to localize the AL subproblem within a compact set $X$ helps in securing \cref{ass:PGoracle:wellposed}.
	Together with properness and lower semicontinuity of $\subLL$, compactness of $X\subset\XX$ ensures that, once localized, the minimization of $\subLL$ is well-posed.
	In the spirit of trust region methods \cite{conn2000trust}, however, it would be necessary to adapt the admissible region $X$ to avoid interfering with the outer AL loop convergence.
	These considerations are left for a future effort.
\end{remark}
	
	\subsection{Descent Methods with an Oracle}\label{sec:PGoracle}
	In this section we detail and discuss an algorithm for the minimization of $\subLL$, which proceeds by taking steps along directions of descent for $\subLL$.
The oracle-equipped descent method is stated in \cref{alg:PGoracle}, which closely patterns \cite[Algorithm 8.1]{birgin2014practical} when variable $z$ is neglected --- recovering a descent method for unconstrained smooth optimization.
We prove well-definedness of the algorithm and investigate its convergence properties in \cref{sec:PGoracle:convergence}.

\begin{algorithm2e}[htb]
	\DontPrintSemicolon
	\KwData{$x^0 \in \XX$, $\varepsilon > 0$}
	\KwResult{$\varepsilon$-$\Upsilon$-stationary point $x^\ast \in \XX$ with certificate $z^\ast \in \FBO(x^\ast)$}
	Select parameters $\alpha, \beta, \theta \in (0,1)$, $\omega, \nu \in (0,1]$ and set $\subLLnm_0 \gets \subLL(x^0)$\;
	\For{$j = 0,1,\ldots$}{
		Select $z^j \in \FBO(x^j)$\label{step:PGoracle:z} and set $\nabla \subLL^j \gets \nabla_x \subLLsmooth(x^j,z^j)$\;
		\lIf{$\| \nabla \subLL^j \| \leq \varepsilon$\label{step:PGoracle:end:if}}{%
			\Return{$(x^\ast, z^\ast) \gets (x^j, z^j)$}
		}%
		Compute $d^j \in \XX$ such that\label{step:PGoracle:d}
		\begin{equation*}
			\innprod{\nabla \subLL^j}{d^j} {}\leq{} - \theta \| \nabla \subLL^j \| \| d^j \|
			\quad\text{and}\quad
			\| d^j \| {}\geq{} \omega \| \nabla \subLL^j \|
		\end{equation*}
		
		Compute $x^{j+1} \in \XX$ and the largest $\gamma_j = \beta^{\ell_j}$, $\ell_j \in \N$, such that\label{step:PGoracle:x}\label{step:PGoracle:gammaLS}
		\begin{equation*}
			\subLL(x^{j+1}) \leq \subLL(x^j + \gamma_j d^j) \leq \subLLnm_j + \alpha \gamma_j \innprod{\nabla \subLL^j}{d^j}
		\end{equation*}

		Set $\subLLnm_{j+1} \gets (1 - \nu) \subLLnm_j + \nu \subLL(x^{j+1})$\label{step:PGoracle:Phi}\;
	}
	\caption{Nonmonotone Descent Method with an Oracle}
	\label{alg:PGoracle}
\end{algorithm2e}%

Some comments are in order.
First, the method requires a starting point $x^0 \in \XX$ for variable $x$ only, not for $z$.
Then, at \cref{step:PGoracle:z}, a surrogate gradient $\nabla \subLL^j \in \XX$ of $\subLL$ at $x^j$ is evaluated: the oracle returns an arbitrary element $z^j \in \FBO(x^j)$ and, in light of the structure \eqref{eq:inner:cost}, we let $\nabla \subLL^j$ be the partial derivative $\partial_x \subLLslack$ at $(x^j,z^j)$.
At \cref{step:PGoracle:d} a suitable search direction $d^j \in \XX$ has to be selected, possibly among many.
The simplest choice is the \emph{steepest} direction given by the negative gradient, i.e. $d^j \coloneqq - \nabla \subLL^j$, which is always valid with respect to the conditions at \cref{step:PGoracle:d}, for any choice of $\theta \in (0,1)$ and $\omega \in (0,1]$.
Another possibility is to consider Newton-like directions, based on local quadratic approximations.
Barzilai-Borwein stepsizes and (limited memory) BFGS directions are popular strategies in this domain; see \cite[Section 5.1.2]{themelis2018forward} for an overview.
Afterwards, the inner loop at \cref{step:PGoracle:gammaLS} seeks an appropriate stepsize $\gamma_j > 0$ such that the update $x^j + \gamma_j d^j$ yields some improvement.
In particular, this backtracking linesearch procedure is herein included as a globalization mechanism.
Finally, with $\nu \in (0,1)$, the update rule at \cref{step:PGoracle:Phi} allows to relax the requirement of a monotone decrease for the objective values, thus reducing conservatism of the linesearch and possibly encouraging longer steps and improved practical performance \cite[Section 8.2]{birgin2014practical}, \cite{themelis2018forward}.
However, a monotone behavior for $\{p(x^j)\}_{j\in\N}$ can be enforced with $\nu \coloneqq 1$.

Crucially, the tentative update and the backtracking condition at \cref{step:PGoracle:gammaLS} for adapting the stepsize $\gamma_k$ are concerned with variable $x$ only.
Also, the oracle evaluation at \cref{step:PGoracle:z} can be formally excluded from the algorithm and incorporated into the \emph{surrogate} gradient evaluation of $\subLL$ at $x^j$.
This is possible thanks to the arbitrariness of the element $z^j$ chosen by the oracle $\FBO$, which allows us to keep variable $z$ implicit.
These observations show that \cref{alg:PGoracle} is equivalent to a descent method applied directly to the minimization of $\subLL$ (with a surrogate gradient), thus making the oracle-equipped algorithm a formal tool only, transparent to the user.
	
	\subsection{Convergence Analysis}\label{sec:PGoracle:convergence}
	In this section we investigate the convergence properties of \cref{alg:PGoracle}, exploiting continuous differentiability of $\subLLsmooth$ and the oracle's optimality to prove sufficient decrease with respect to $x$.
We are going to show first that \cref{alg:PGoracle} is well defined, namely that each iteration terminates in finite time.
This is an essential theoretical and practical condition that an implementable algorithm must satisfy.
\begin{mybox}
	\begin{lemma}[Well definedness]\label{thm:PGoracle:finite}%
		Let \cref{ass:P} hold and consider the iterates generated by \cref{alg:PGoracle}.
		Then,
		\begin{enumerate}[label=(\roman*)]
			\item\label{thm:PGoracle:finite:LS}
			at every iteration, the number of backtrackings at \cref{step:PGoracle:gammaLS} is finite;
			\item\label{thm:PGoracle:finite:descent}
			at the end of the $j$-th iteration, one has
			\begin{equation*}
				\subLL(x^{j+1})
				\leq
				\Phi_j + \delta_j
				\qquad\text{and}\qquad
				\subLL(x^{j+1})
				\leq
				\Phi_{j+1}
				\leq
				\Phi_j + \nu \delta_j
			\end{equation*}
			where $\delta_j \coloneqq \alpha \gamma_j \innprod{\nabla \subLL^j}{d^j} < 0$;
			\item\label{thm:PGoracle:finite:sublevel}
			every iterate $x^j$ remains within the sublevel set $\{ x\in\XX \,\vert\, \subLL(x) \leq \Phi_0 < \infty \}$.
		\end{enumerate}
	\end{lemma}
\end{mybox}
\begin{proof}
	Recall the basic properties collected in \cref{thm:PGoracle:properties}.
	\begin{itemize}
		\item[\ref{thm:PGoracle:finite:LS}]
		Let us discard the trivial case by assuming that $\| \nabla \subLL^j \| > \varepsilon$.
		Since $\theta \in (0,1)$ and $\omega \in (0,1]$, the choice $d^j \gets - \nabla \subLL^j$ is always viable, hence \cref{step:PGoracle:d} is well-defined.
		Denote $x^{j+} \coloneqq x^j + \gamma_j d^j$ and $z^{j+} \in \FBO(x^{j+})$.
		Then,
		\begin{equation*}
			\subLL(x^{j+})
			=
			\subLLsmooth(x^{j+}, z^{j+}) + g(z^{j+})
			\leq
			\subLLsmooth(x^{j+}, z) + g(z)
		\end{equation*}
		for all $z\in\YY$.
		In particular, $\subLL(x^{j+}) \leq \subLLsmooth(x^{j+}, z^j) + g(z^j)$.
		By continuous differentiability of $\subLLsmooth$, combining \cref{step:PGoracle:d,step:PGoracle:z} yields
		\begin{equation*}
			\lim\limits_{\gamma\to 0} \frac{\subLLsmooth(x^j + \gamma d^j,z^j) - \subLLsmooth(x^j,z^j)}{\gamma}
			=
			\innprod{ \nabla_x  \subLLsmooth(x^j,z^j) }{ d^j }
			=
			\innprod{ \nabla \subLL^j }{ d^j }
			<
			0 .
		\end{equation*}
		Thus,
		\begin{equation*}
			\lim\limits_{\gamma\to 0} \frac{\subLLsmooth(x^j + \gamma d^j,z^j) - \subLLsmooth(x^j,z^j)}{\gamma \innprod{ \nabla \subLL^j }{ d^j }}
			=
			1 .
		\end{equation*}
		Therefore, since $\alpha \in (0,1)$, for some $\gamma_j = \beta^{\ell_j} > 0$ small enough we have that
		\begin{equation*}
			\frac{\subLLsmooth(x^j + \gamma_j d^j,z^j) - \subLLsmooth(x^j,z^j)}{\gamma_j \innprod{ \nabla \subLL^j }{ d^j }}
			\geq
			\alpha .
		\end{equation*}
		As $\innprod{ \nabla \subLL^j }{ d^j } < 0$ by \cref{step:PGoracle:d} and $\gamma_j > 0$, we deduce that
		\begin{equation*}
			\subLLsmooth(x^{j+}, z^j)
			\leq
			\subLLsmooth(x^j, z^j) + \alpha \gamma_j \innprod{ \nabla \subLL^j }{ d^j }
			.
		\end{equation*}
		Now, invoking the oracle to get $z^{j+} \in \FBO(x^{j+})$ we obtain
		\begin{multline*}
			\subLL(x^{j+})
			=
			\subLLsmooth(x^{j+}, z^{j+}) + g(z^{j+})
			\leq
			\subLLsmooth(x^{j+}, z^j) + g(z^j) \\
			\leq
			\subLLsmooth(x^j, z^j) + g(z^j) + \alpha \gamma_j \innprod{ \nabla \subLL^j }{ d^j }
			=
			\subLL(x^j) + \alpha \gamma_j \innprod{ \nabla \subLL^j }{ d^j }
			.
		\end{multline*}
		Then, finite termination of the inner loop follows from the bound $\subLLnm_j \geq \subLL(x^j)$ valid for all $j\in\N$.
		The latter is obtained by induction: if $j=0$, it is $\subLLnm_j = \subLL(x^j)$ by construction; if $j>0$, the bound is due to \cref{step:PGoracle:Phi}.
		\item[\ref{thm:PGoracle:finite:descent}]
		The first inequality follows by combining \ref{thm:PGoracle:finite:LS} with the failure of the condition at \cref{step:PGoracle:gammaLS}.
		Together with \cref{step:PGoracle:Phi}, this gives the lower bound
		\begin{multline*}
			\subLLnm_{j+1}
			=
			(1 - \nu) \subLLnm_j + \nu \subLL(x^{j+1})
			\geq
			(1 - \nu) [\subLL(x^{j+1}) - \delta_j] + \nu \subLL(x^{j+1})\\
			=
			\subLL(x^{j+1}) - (1 - \nu) \delta_j
			\geq
			\subLL(x^{j+1})
		\end{multline*}
		as well as the upper bound
		\begin{equation*}
			\subLLnm_{j+1}
			=
			(1 - \nu) \subLLnm_j + \nu \subLL(x^{j+1})
			\leq
			(1 - \nu) \subLLnm_j + \nu (\subLLnm_j + \delta_j)
			=
			\subLLnm_j + \nu \delta_j .
		\end{equation*}
		\item[\ref{thm:PGoracle:finite:sublevel}]
		That $\Phi_0 = \subLL(x^0) < \infty$ follows from the fact that $z^0 \in \FBO(x^0) \in \dom g \ne \emptyset$.
		Then, the assertion is a direct consequence of \ref{thm:PGoracle:finite:descent}, by telescoping $\subLL(x^j) \leq \subLLnm_j \leq \ldots \leq \Phi_0$.
		\qedhere
	\end{itemize}
\end{proof}

We next consider an asymptotic analysis of the algorithm.
\begin{mybox}
	\begin{lemma}\label{thm:PGoracle:asymp:cost}%
		Let \cref{ass:P,ass:PGoracle} hold and consider a sequence of iterates $\{x^j\}_{j\in\N}$ generated by \cref{alg:PGoracle} with $\varepsilon = 0$.
		Then, the sequences $\{ \Phi_j \}_{j\in\N}$ and $\{\subLL(x^j)\}_{j\in\N}$ converge to a finite value $\subLL_\ast \geq \inf \subLL$, the former from above.
	\end{lemma}
\end{mybox}
\begin{proof}
	Regarding $\{\Phi_j\}_{j\in\N}$, the assertion follows from \cref{thm:PGoracle:finite}\ref{thm:PGoracle:finite:descent} and \cref{ass:PGoracle}.
	Then, rearranging the update rule at \cref{step:PGoracle:Phi}, one can observe that $\subLL(x^{j+1}) = \nu^{-1} (\Phi_{j+1} - \Phi_j) + \Phi_j$ holds. Taking the limit for $j \to \infty$, by $\nu > 0$ the convergence of $\{\Phi_j\}_{j\in\N}$ implies that of $\{\subLL(x^j)\}_{j\in\N}$ to the same finite value $\subLL_\ast$.
	\qedhere
\end{proof}

The following \cref{thm:PGoracle:convergenceStationarity} provides a guarantee of global convergence to $\Upsilon$-stationary points, namely that any accumulation point of the iterates generated by \cref{alg:PGoracle} is $\Upsilon$-stationary relative to the minimization of $\subLL$.
\begin{mybox}
	\begin{theorem}\label{thm:PGoracle:convergenceStationarity}
		Let \cref{ass:P,ass:PGoracle} hold and consider a sequence of iterates $\{x^j\}_{j\in\N}$ generated by \cref{alg:PGoracle} with $\varepsilon = 0$.
		Let $x^\ast \in \XX$ be an accumulation point of $\{x^j\}_{j\in\N}$.
		Then, $x^\ast$ is $\Upsilon$-stationary for $\subLL$, namely there exists $z^\ast \in \FBO(x^\ast)$ such that $\nabla_x \subLLsmooth(x^\ast, z^\ast) = 0$.
	\end{theorem}
\end{mybox}
\begin{proof}
	If \cref{alg:PGoracle} terminates after finitely many iterations returning $(x^j,z^j)$, then it must be that $\| \nabla p^j \| = 0$, hence $\nabla_x \subLLsmooth(x^j,z^j) = 0$.
	Since $z^j \in \FBO(x^j)$ for all $j\in\N$, then $x^j$ would be $\Upsilon$-stationary.

	Now suppose instead that \cref{alg:PGoracle} takes infinitely many iterations.
	Let $\{x^j\}_{j\in J}$, $J \coloneqq \{j_0,j_1,j_2\ldots\} \subseteq \N$, be a subsequence such that $x^j \to_J x^\ast$.
	Then, continuity of $\subLL$ by \cref{thm:PGoracle:properties}\ref{ass:PGoracle:subLL} gives $\subLL(x^j) \to_J \subLL(x^\ast)$.
	Owing to the linesearch condition at \cref{step:PGoracle:gammaLS}, since $j_{i+1} \geq j_i + 1$, we have that
	\begin{equation*}
		\subLL(x^{j_{i+1}})
		\leq
		\subLL(x^{j_i + 1})
		\leq
		\subLLnm_{j_i} + \alpha \gamma_{j_i} \innprod{\nabla \subLL^{j_i}}{d^{j_i}}
		<
		\subLLnm_{j_i}
		.
	\end{equation*}
	Then, \cref{thm:PGoracle:asymp:cost} implies that $\gamma_j \innprod{\nabla \subLL^j}{d^j} \to_J 0$, since $\alpha > 0$, and, due to \cref{step:PGoracle:d}, it must be that $\gamma_j \| \nabla \subLL^j \| \| d^j \| \to_J 0$.
	As any term in the limit could vanish, there are three (non exclusive) cases.
	Now we argue that, up to extracting a subsequence, it is $\nabla \subLL^j \to_J 0$ in every case.
	First, the claim is trivial if $\nabla \subLL^j$ vanishes.
	Second, if $d^j \to_J 0$, conditions at \cref{step:PGoracle:d} imply that $\nabla \subLL^j$ vanishes too, giving the claim.
	Third, if $\gamma_j \to_J 0$, we proceed by showing that there is a longer stepsize $\hat{\gamma}_j$ that also tends to zero, along which the function did not decrease enough, and that this is possible only if the gradient vanishes.

	Let us start by noticing the failure of attempts at \cref{step:PGoracle:gammaLS}.
	The previous (hence rejected) stepsize $\hat{\gamma}_j \coloneqq \beta^{-1} \gamma_j = \beta^{\ell_j - 1}$ gives
	\begin{equation}
		\subLL(x^j + \hat{\gamma}_j d^j)
		{}>{}
		\subLLnm_j + \alpha \hat{\gamma}_i \innprod{\nabla \subLL^j}{d^j}
		{}\geq{}
		\subLL(x^j) + \alpha \hat{\gamma}_j \innprod{\nabla \subLL^j}{d^j} , \label{eq:PG:failed_stepsize}
	\end{equation}
	where the last inequality follows from the bound $\subLLnm_j \geq \subLL(x^j)$ of \cref{thm:PGoracle:finite}\ref{thm:PGoracle:finite:descent}.
	In the case $\gamma_j \to_J 0$, it must be $\hat{\gamma}_j \to_J 0$ too, since $\ell_j \to _J 0$.
	Rearranging terms and using the structure of $\subLLslack$ in \eqref{eq:inner:cost}, the minimizing property of any $z^j \in \FBO(x^j)$ yields
	\begin{multline*}
		\alpha \innprod{\nabla \subLL^j}{d^j}
		<
		\frac{\subLL(x^j + \frac{\gamma_j}{\beta} d^j) - \subLL(x^j)}{\frac{\gamma_j}{\beta}}\\
		\leq
		\frac{\subLLsmooth(x^j + \frac{\gamma_j}{\beta} d^j, z^j) + g(z^j) - \subLLsmooth(x^j,z^j) - g(z^j)}{\frac{\gamma_j}{\beta}}\\
		=
		\frac{\subLLsmooth(x^j + \frac{\gamma_j}{\beta} d^j, z^j) - \subLLsmooth(x^j,z^j)}{\frac{\gamma_j}{\beta}}
		.
	\end{multline*}
	Defining $s^j \coloneqq \hat{\gamma}_j d^j$ for all $j\in\N$, we have that $\|s^j\| \to _J 0$.
	Possibly extracting a subsequence and redefining $J$, let $s^\ast \in \XX$ be such that $s^j/\|s^j\| \to_J s^\ast$.
	Now, for all $j\in J$, by the mean value theorem there exists $\zeta_j \in [0,1]$ such that
	\begin{align*}
		\innprod{ \nabla_x \subLLsmooth(x^j + \zeta_j s^j, z^j) }{s^j}
		{}={}&
		\subLLsmooth(x^j + s^j, z^j) - \subLLsmooth(x^j, z^j) \\
		{}={}&
		\subLLsmooth(x^j + s^j, z^j) + g(z^j) - \subLLsmooth(x^j, z^j) - g(z^j) \\
		{}\geq{}&
		\subLL(x^j + s^j) - \subLL(x^j) \\
		{}>{}&
		\alpha \innprod{\nabla \subLL^j}{s^j}
		{}={}
		\alpha \innprod{\nabla_x \subLLsmooth(x^j, z^j)}{s^j} ,
	\end{align*}
	where the inequalities are due to the oracle's optimality and \eqref{eq:PG:failed_stepsize}.
	Dividing by $\|s^j\|$ and rearranging give
	\begin{equation*}
		\frac{ \innprod{ \nabla_x \subLLsmooth(x^j + \zeta_j s^j, z^j) - \alpha \nabla_x \subLLsmooth(x^j, z^j) }{s^j} }{\|s^j\|}
		>
		0 .
	\end{equation*}
	Owing to $\zeta_j s^j \to_J 0$ and continuity of $\nabla \subLLsmooth$, taking the limit yields
	\begin{multline*}
		\limsup\limits_{j \to_J \infty} \frac{ \innprod{ \nabla_x \subLLsmooth(x^j + \zeta_j s^j, z^j) - \alpha \nabla_x \subLLsmooth(x^j, z^j) }{s^j} }{\|s^j\|} \\
		{}={}
		\limsup\limits_{j \to_J \infty} (1 - \alpha) \frac{ \innprod{ \nabla_x \subLLsmooth(x^j, z^j) }{s^j} }{\|s^j\|}
		{}\leq{} 0 ,
	\end{multline*}
	where the inequality is due to $\innprod{ \nabla_x \subLLsmooth(x^j, z^j) }{s^j} = \hat{\gamma}_j \innprod{ \nabla \subLL^j }{ d^j } < 0$ for all $j\in\N$, by \cref{step:PGoracle:d}, and $\alpha < 1$.
	Recalling that $s^j/\|s^j\| \to_J s^\ast$, the inequalities above demonstrate that it is necessarily $\innprod{ \nabla_x \subLLsmooth(x^j, z^j) }{s^\ast} = \innprod{ \nabla \subLL^j }{ s^\ast } \to_J 0$.
	Then, by the conditions at \cref{step:PGoracle:d}, it must be $\nabla \subLL^j \to_J 0$.

	Finally, owing to $\nabla \subLL^j \coloneqq \nabla_x \subLLsmooth(x^j, z^j) \to_J 0$ in all cases, $z^j \in \FBO(x^j)$ for all $j\in J$, $x^j \to_J x^\ast$, and continuity of $\nabla \subLLsmooth$, we conclude that $\nabla_x \subLLsmooth(x^\ast,z^\ast) = 0$ for some $z^\ast \in \FBO(x^\ast)$, proving the claim that $x^\ast$ is $\Upsilon$-stationary for $\subLL$.
	\qedhere
\end{proof}

	\subsection{Termination Criteria}\label{sec:PGoracle:termination}
	\cref{alg:PGoracle} would run indefinitely without \cref{step:PGoracle:end:if}, generating an infinite sequence of iterates.
The merit and objective values tend to decrease along the iterations and settle at some finite value, as established by \cref{thm:PGoracle:finite}\ref{thm:PGoracle:finite:descent} and \cref{thm:PGoracle:asymp:cost}, respectively.
Suitable termination criteria need to be inserted for stopping and returning an iterate that, in some sense, is satisfactory relatively to the minimization of $\subLL$.

Since $z^j \in \FBO( x^j )$ for all $j \in \N$, it suffices to check the (approximate) stationarity of $x^j$.
Thus, we equip \cref{alg:PGoracle} with the termination condition $\left\| \nabla_x \subLLsmooth(x^j, z^j) \right\| \leq \varepsilon$, which encodes $\varepsilon$-$\Upsilon$-stationarity of $x^j$ relative to $\subLL$, by the specific structure of $\subLLslack$ in \eqref{eq:inner:cost}.
Thanks to \cref{thm:PGoracle:convergenceStationarity} we can now easily infer termination of \cref{alg:PGoracle} in finitely many steps for any tolerance $\varepsilon > 0$, which confirms that its output complies with the requirements for the outer AL framework at \cref{step:ALM:subproblem} of \cref{alg:ALM}.
\begin{mybox}
	\begin{corollary}\label{thm:PG:finitetermination}
		Let \cref{ass:P} hold and consider a sequence of iterates $\{x^j,z^j\}_{j\in\N}$ generated by \cref{alg:PGoracle} with $\varepsilon > 0$.
		Suppose that the sequence $\{x^j\}_{j\in\N}$ admits an accumulation point.
		Then, in finitely many steps an $\varepsilon$-$\Upsilon$-stationary point $x^\ast \in \dom\subLL$ is returned with certificate $z^\ast \in \FBO(x^\ast)$.
	\end{corollary}
\end{mybox}
\begin{proof}
	That the algorithm terminates in finitely many iterates, say $j_\ast$ many, follows from \cref{thm:PGoracle:convergenceStationarity}, since $\| \nabla_x \subLLsmooth(x^j,z^j) \| = \| \nabla \subLL^j \| \to_J 0$ along any convergent subsequence $\{x^j\}_{j\in J}$.
	Since $z^j \in \FBO(x^j)$ for all $j\in\N$, it follows that the output $x^\ast \coloneqq x^{j_\ast}$ with certificate $z^\ast \coloneqq z^{j_\ast}$ satisfies $z^\ast \in \FBO(x^\ast)$.
	Since the magnitude of $\nabla_x \subLLsmooth(x^\ast,z^\ast)$ is no more than $\varepsilon$ by the termination criterion, the pair $(x^\ast,z^\ast)$ is $\varepsilon$-$\Upsilon$-stationary for $\subLL$.
	\qedhere
\end{proof}
Note that the existence of limit points is assumed, and not guaranteed, in \cref{thm:PG:finitetermination}.
This can be satisfied by employing, in general, (level) boundedness or coercivity arguments.

\section{Numerical Example and Comparison}\label{sec:numresults}
This section showcases the flexibility offered by the (nonconvex) generalized programming framework and presents some numerical results that support the implicit approach developed in \cref{sec:ALM}.
We consider as illustrative example a well-known academic problem with vanishing constraints.

We have implemented the proposed approach in the Augmented Lagrangian Proximal Solver (\alps{}), as part of the open-source software package \href{https://github.com/aldma/Bazinga.jl}{Bazinga.jl} \cite{demarchi2020bazinga}, which includes the example discussed below.
As our goal is to provide only a motivation and numerical support for the implicit approach, we do not include implementation details.
Nonetheless, we shall point out that it closely patterns that of \als{} \cite{demarchi2023constrained}, an AL solver based on an explicit treatment of slack variables.
In particular, the AL subproblems at \cref{step:ALM:subproblem} of \cref{alg:ALM} are solved by default using the accelerated proximal gradient solver PANOC$^+$ \cite{demarchi2022proximal}.
This choice is convenient in our numerical tests to obtain comparable results across different approaches and formulations.
As a proximal gradient method, PANOC$^+$ can handle the two-term structure of the explicit AL function; it is in fact the default inner solver in \als{} \cite{demarchi2023constrained}.
However, PANOC$^+$ can act as a gradient-based descent method too, and thus can be adopted as inner solver within \cref{alg:ALM}, in view of \cref{sec:PGoracle}.

	\subsection{Problem with Vanishing Constraints}\label{sec:numresults:academicmpvc}
\begin{figure}[htb]
	\centering
	\includegraphics{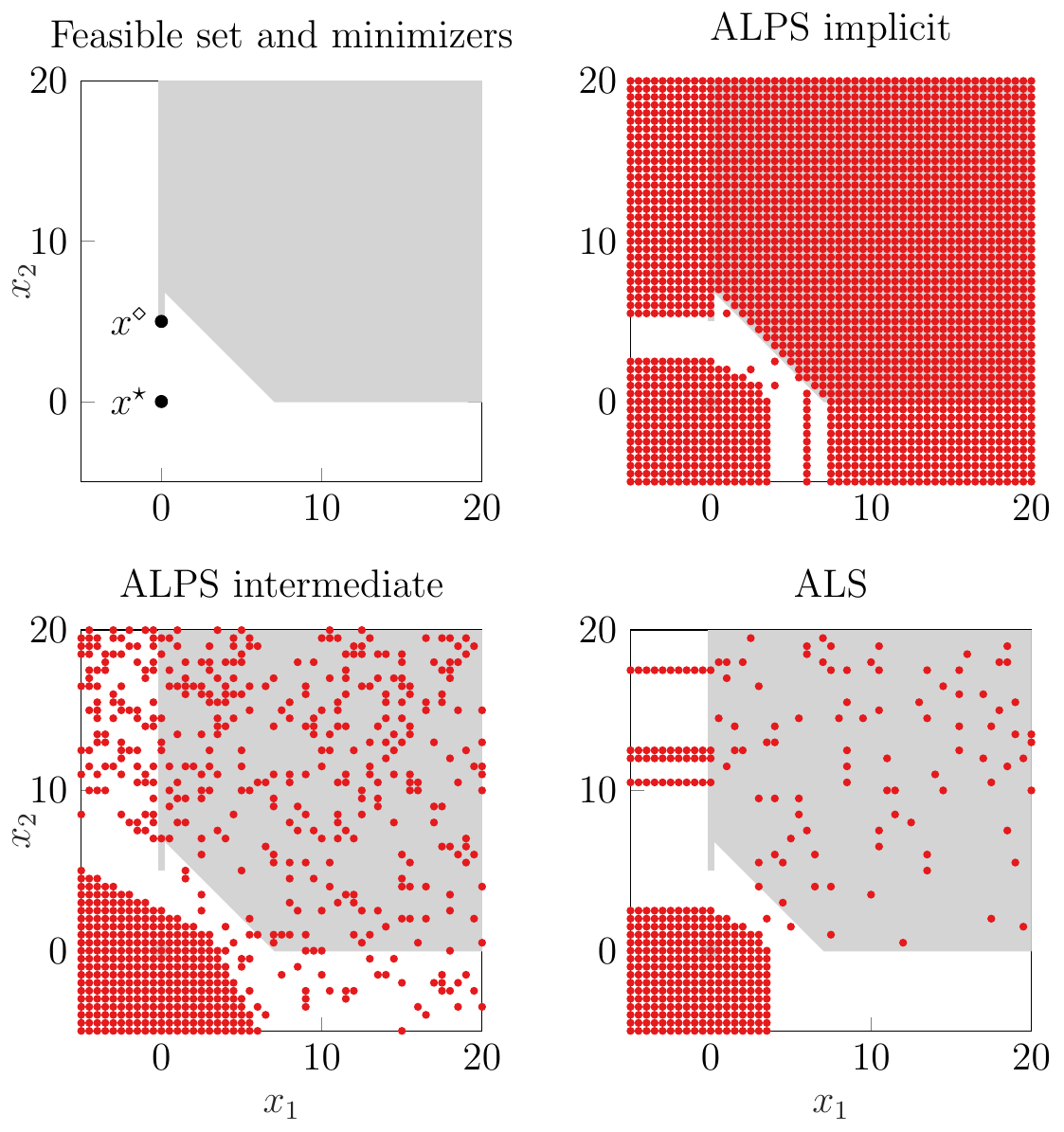}
	\caption{Comparison of solutions obtained using different methods and formulations. Feasible set (gray background) and grid of (2601 uniformly spaced) starting points and corresponding termination mark depending on the returned solution: global minimizer $x^\star = (0,0)$ (red circle) or local minimizer $x^\diamond = (0,5)$ (omitted mark).}
	\label{fig:mpvca}
\end{figure}

	Let us consider the following two-dimensional problem, which arises from truss topology optimization and is a well-known academic problem; it has been thoroughly discussed in \cite[Section 9.5.1]{hoheisel2009mathematical}.
\begin{align}
	\minimize_{x \in \R^2}
	\quad{}{}&
	4 x_1 + 2 x_2 \label{eq:mpvca}\\
	\stt
	\quad{}{}&
	x_1 \geq 0 ,\quad
	x_2 \geq 0 \nonumber \\
	{}{}&
	x_1 > 0 {}\quad\Rightarrow\quad{} x_1 + x_2 \geq 5 \sqrt{2} \nonumber \\
	{}{}&
	x_2 > 0 {}\quad\Rightarrow\quad{} x_1 + x_2 \geq 5 \nonumber
\end{align}
Due to their structure, the logical implications give rise to the so called vanishing constraints, usually reformulated by requiring
$x_1 (x_1 + x_2 - 5 \sqrt{2}) {}\geq{} 0$
and
$x_2 (x_1 + x_2 - 5) {}\geq{} 0$
along with $x_1, x_2 \geq 0$.
The feasible set of \eqref{eq:mpvca} consists of the union of an unbounded polyhedron $\{ (x_1,x_2) \in \R_+^2 \,\vert\, x_1+x_2\geq5\sqrt{2}\}$, of an attached line segment $\{ (0, x_2 ) \in \R^2 \,\vert\, 5 \leq x_2 \leq 5\sqrt{2} \}$, and of the isolated point $(0,0)$.
It admits a global minimum at $x^\star = (0,0)$ and a local minimum at $x^\diamond = (0,5)$; see \cref{fig:mpvca}.

There exist several ways to cast \eqref{eq:mpvca} into the form of \eqref{eq:P}.
We are going to consider three different formulations, with increasing degree of explicitness for treating constraints: implicit, intermediate (or partially explicit, avoiding slack variable components for trivial constraints), and (fully) explicit.
Let us define the cost function $f(x) \coloneqq 4 x_1 + 2 x_2$.
Moreover, let us conveniently denote by $\Dvc \subset \R^2$ the set associated to a vanishing constraint: $\Dvc \coloneqq \{ (a,b) \in \R^2 \,\vert\, a \geq 0 , a b \geq 0 \}$; it is nonempty, closed, nonconvex, and its (set-valued) projection $\proj_D$ can be easily evaluated.
\begin{itemize}
	\item An \emph{implicit} representation of \eqref{eq:mpvca} in the form \eqref{eq:P} is given by the data functions
	\begin{align*}
		c_{\text{impl}}(x)
		{}\coloneqq{}&
		(x_1 ,\, x_1 + x_2 - 5 \sqrt{2} ,\, x_2 ,\, x_1 + x_2 - 5)^\top \\
		g_{\text{impl}}(c)
		{}\coloneqq{}&
		\indicator_{\Dvc}(c_1,c_2) + \indicator_{\Dvc}(c_3,c_4)
		.
	\end{align*}
	With size $(n,m)=(2,4)$, this seems the most compact formulation for \eqref{eq:mpvca}.
	\item An \emph{explicit} reformulation of \eqref{eq:mpvca} as \eqref{eq:P} is obtained by introducing slack variables as in \eqref{eq:ALMz:P}.
	This leads to a problem with size $(n,m)=(6,8)$.
	Notice that this explicit reformulation is carried out internally by \als{} when fed with the implicit formulation, as a lifted representation of the original problem \cite[Section 3]{demarchi2023constrained}.
	\item An \emph{intermediate} representation of \eqref{eq:mpvca} into the form of \eqref{eq:P} is found by avoiding slack variable components for trivial constraints.
	This intermediate formulation has size $(n,m)=(4,6)$, with data functions
	\begin{align*}
		c_{\text{inte}}(x)
		{}\coloneqq{}&
		(x_1 ,\, x_1 + x_2 - 5 \sqrt{2} - x_3 ,\, x_2 ,\, x_1 + x_2 - 5 - x_4 ,\, x_3 ,\, x_4)^\top \\
		g_{\text{inte}}(c)
		{}\coloneqq{}&
		\indicator_{\Dvc}(c_1,c_5) + \indicator_{\Dvc}(c_3,c_6) + \indicator_0(c_2) + \indicator_0(c_4)
		.
	\end{align*}
\end{itemize}

We execute \alps{} on all problem representations.
The outcome of \als{} \cite{demarchi2023constrained} on the implicit formulation is also recorded, for comparison with the explicit one on \alps{}.
We consider a grid of $51^2 = 2601$ starting points $x^0$ in $[-5,20]^2$ and let the initial Lagrange multiplier estimate be $y^0 = 0$; slack variables are initialized so as to satisfy the linear constraints arising from the reformulation.
Both solvers return a successful status for all instances, always indicating either the global or local minimizer (within $10^{-6}$ in Euclidean distance).
In \cref{fig:mpvca} we depict the feasible region of \eqref{eq:mpvca} and display these results in more detail.
Each starting point $x^0$ is given a mark, indicating the termination point which has been reached by starting from $x^0$.
Since the set $\Dvc$ is nonconvex, the iterates of \cref{alg:ALM} depend on the arbitrary choices made when the projection onto $\Dvc$ is not a singleton.
By selecting different projections, we observed analogous numerical results, albeit not identical.

\cref{tab:mpvca} collects statistics that summarize the numerical results.
Considering the \emph{implicit} formulation, \alps{} returns the global minimizer $x^\star = (0,0)$ in most cases ($93\%$).
With the \emph{intermediate} formulation, \alps{} seems to consistently find $x^\star$ only when starting nearby ($28\%$), and even more so with the fully \emph{explicit} formulation ($19\%$) or with \als{} ($16\%$).
Despite the fact that solutions and stationary points remain unaffected by the lifting \cite[Lemma 3.1]{demarchi2023constrained}, there is a significant gap in terms of objective values attained in practice.
But the advantages could go beyond this aspect, as the implicit approach seems extremely effective in reducing the (median) number of iterations and runtime too.
We are drawn to suspect that this is the oracle's doing.

\begin{table}[htb]
	\centering
	\caption{Comparison of solutions obtained and performances using different methods and formulations. Statistics on the number of (outer and cumulative inner) iterations, elapsed runtime, and how often the global minimizer $x^\star$ is returned.}
	\label{tab:mpvca}
	\medskip
	\begin{tabular}{cc|cccc}
		& & {\alps{}} & {\alps{}} & {\alps{}} & {\als{}} \\
		& & \emph{implicit} & \emph{intermediate} & \emph{explicit} & \\
		\hline
		outer & median & 10 & 10 & 10 & 10 \\
		iterations & max & 10 & 10 & 10 & 10 \\
		\hline
		total inner & median & 10 & 77 & 117 & 117 \\
		iterations & max & 802 & 598 & 217 & 205 \\
		\hline
		runtime & median & 0.17 & 2.1 & 1.9& 3.0 \\
		{[ms]} & max & 155 & 130 & 46 & 36 \\
		\hline
		global & abs & 2413 & 721 & 496 & 405 \\
		minimizer & rel & 92.77\% & 27.72\% & 19.07\% & 15.57\% \\
	\end{tabular}
\end{table}

\section{Concluding Remarks}\label{sec:conclusions}
The theoretical investigation of generalized optimization in the fully nonconvex setting has led to a significant extension of the augmented Lagrangian (AL) framework.
With the help of a lifted reformulation, KKT-type stationarity conditions have been derived in Lagrangian terms and complemented with asymptotic and regularity notions.
The implicit AL function and subproblems have been defined, marginalizing their explicit counterparts with respect to slack variables, seemingly sacrificing some regularity due to set-valuedness.
Then, we devised the concept of $\Upsilon$-stationarity, rooted on a parametric interpretation, for characterizing the subproblems' solutions.
Although weaker than M-stationarity for the implicit AL function, $\Upsilon$-stationarity proved to be a strong enough tool for investigating asymptotic properties. 
We considered a safeguarded AL algorithm and established opportune convergence guarantees, on par with those in nonlinear programming.

Embracing an implicit philosophy was not only a theoretical curiosity, though.
We developed and characterized a solver for the implicit AL subproblems, which is indeed capable of generating approximate $\Upsilon$-stationary iterates.
Equipped with a proximal oracle, it handles fewer decision variables, maintaining slack variables under the hood.
This appears to produce a tighter model and catalyze improved practical performance, compared to more explicit approaches, as witness by numerical results on an illustrative problem.

Several issues and questions remain open.
It is a common practice to assume that the AL subproblem are well-posed, possibly based on coercivity or boundedness arguments.
In general, however, to overcome this concern one could investigate the combination of AL schemes with trust-region strategies, as hinted at in \cref{rem:abstractConstraints}.
Another topic of interest is that of attentive convergence, whose role should be closely inspected.
It remains unclear whether it can be discarded or not, in relation to asymptotic stationarity conditions.
Finally, future works should elaborate on a characterization of stationary points as saddle points of the implicit AL function, a local duality theory for generalized optimization, and the relationships with proximal point algorithms.

\ifpreprint
	\subsection*{Acknowledgements}
	\TheAcknowledgementsADM
\fi

{%
\ifpreprint
	\phantomsection
	\addcontentsline{toc}{section}{References}%
	\small
	\bibliographystyle{plain}
\fi
\bibliography{biblio}
}%

\end{document}